\documentclass[preprint,11pt]{article}

\usepackage{type1cm} 
\usepackage{graphicx}   
\usepackage[bottom]{footmisc}   
\usepackage{newtxtext}        
\usepackage{newtxmath}   
\usepackage{amsmath,amsfonts}
\usepackage{titlesec}
\usepackage{booktabs}
\usepackage{multirow}
\usepackage{bbm}
\usepackage{algorithm}
\usepackage{algorithmic} 
\usepackage{todonotes}
\usepackage{hyperref}
\usepackage{authblk}
\usepackage[a4paper, total={7in, 9in}]{geometry}

\def\R{{\mathbb R}}

\newcommand{\ko}{k_o}
\newcommand{\kd}{k_d}
\newtheorem{prop}{Proposition}

\begin{document}

\title{Data/moment-driven approaches for fast predictive control\\ of collective dynamics}

\author[1]{Giacomo Albi}
\author[2]{Sara Bicego}
\author[3]{Michael Herty}
\author[2]{Yuyang Huang}
\author[2]{Dante Kalise}
\author[3]{Chiara Segala}

\affil[1]{ Albi Dipartimento di Informatica, Universit\`a di Verona, Verona, Italy.
	\texttt{e-mail:s.giacomo.albi@univr.it}}
\affil[2]{Department of Mathematics, Imperial College London, United Kingdom. \texttt{e-mail:s.bicego21,yuyang.huang21,dkaliseb@imperial.ac.uk}}
\affil[3]{IGPM, RWTH Aachen University, Templergraben 55, D-52062 Aachen, Germany. \texttt{e-mail:herty,segala@igpm.rwth-aachen.de}}

\maketitle
	
\let\thefootnote\relax\footnotetext{
M.H. and C.S. thank the Deutsche Forschungsgemeinschaft (DFG, German Research Foundation) for the financial support 
442047500/SFB1481 within the projects B04, B05, B06, and SPP 2298 Theoretical Foundations of Deep Learning  within the Project(s) HE5386/23-1, Meanfield Theorie zur Analysis von Deep Learning Methoden (462234017). G.A. is member of GNCS of INdAM, and thanks the support of MIUR-PRIN Project 2022 PNRR No. P2022JC95T, Data-driven discovery and control of multi-scale interacting artificial agent systems. This research was supported by the UK Engineering and Physical Sciences Research Council (EPSRC) grant EP/T024429/1
}
	
\abstract{
Feedback control synthesis for large-scale particle systems is reviewed in the framework of model predictive control (MPC). The high-dimensional character of collective dynamics hampers the performance of traditional MPC algorithms based on fast online dynamic optimization at every time step. Two alternatives to MPC are proposed. First, the use of supervised learning techniques for the offline approximation of optimal feedback laws is discussed. Then, a procedure based on sequential linearization of the dynamics based on macroscopic quantities of the particle ensemble is reviewed. Both approaches circumvent the online solution of optimal control problems enabling fast, real-time, feedback synthesis for large-scale particle systems. Numerical experiments assess the performance of the proposed algorithms.
}

\section{Introduction}\label{sec:introduction}

The description of collective phenomena by means of agent-based systems is a consolidate line of research in several scientific fields such as socio-economics, biology, and robotics, as highlighted in \cite{bellomo20review, MR2974186, MR2165531, MR3119732, MR2861587, MR2580958,Giselle,BioKing}. 
The mathematical modeling of agent-based systems consists with a large ensemble of differential equations describing the evolution of the trajectories of $N$ particles, whose dynamics is ruled by the superposition of binary interaction forces. These interactions encode nonlinear mechanisms based on different social principles like attraction, repulsion, and alignment. One distinctive aspect of these models lies in their capability of reproducing intricate dynamics such as self-organized patterns, including consensus, flocking, and milling, as discussed in \cite{MR2887663, MR2247927,cucker2007emergent,d2006self,motsch2014heterophilious}.

A question of particular interest pertains to determining external control actions capable of influencing the collective behavior of these systems, which can be properly formulated as an optimal control problem. This understanding is crucial for various practical applications, facilitating the development of tailored strategies such as collision-avoidance protocols for swarm robotics \cite{CKPP19,KPAsurvey15,MR3157726,Meurer}, pedestrian evacuation strategies for crowd dynamics \cite{MR3542027,MR3308728,dyer2009leadership,MR4046175}, interventions assessment in traffic management \cite{MR3948232,han2017resolving,stern2018dissipation}, and learning of opinion dynamics \cite{albi2023data,MR3268062,Garnier}
However, the high dimensionality given by the state space and the large number of agents does not usually allow the synthesis of optimal feedback control by means of direct methods, see for example \cite{CFPT15,Bailo_2018}.
Hence, the development of novel efficient algorithms are key to handling real-time applications, which often stumble upon cumbersome software and hardware requirements.

Our aim is to approach this task from various perspectives by reviewing techniques developed in the context of Model Predictive Control (MPC) methods for approximating control strategies for agent-based systems. One of the key challenges associated with MPC strategies is its computational complexity, given the need for subsequent solutions to several optimal control problems. This can be computationally intensive, especially for high-dimensional systems and long prediction horizons. 

The research community has been introducing promising remedies for these problems, involving the simplification of either the predictive models or the optimization problems, to reduce the computational complexity while maintaining acceptable control performance. Among the various proposed methodologies, we cite Reduced-Order MPC \cite{RoM1,RoM2,RoM3}, linearization techniques \cite{LinMPC}, lookup tables (e.g. explicit MPC \cite{ExplicitMPC,ExplicitMPC1}), fast-update schemes \cite{FastMPC}, and approximation techniques for the optimization tasks, mostly based on machine learning algorithms. For instance, Learning-based MPC addresses the data-driven approximation of possibly different elements of the MPC formulation, such that the overall performance is improved.  For these methods, a fundamental categorisation can be made between offline and online paradigms. In offline learning, the model targets the controller design and is trained using synthetic samples, or data obtained from previous simulations \cite{offline1,offline2}. 
In online learning instead, the model is trained to approximate the system dynamics, and data are collected while the system evolves, hence the model is actively adjusted during the closed-loop operation. Among other methods in this class, we find Gaussian Processes \cite{GP1,GP2,Gp3} -- providing a non-parametric and probabilistic model of the system's behavior -- and Reinforcement Learning \cite{RL1,RL2,RL3}, particularly useful when the dynamics are complex, unknown, or hard to model accurately using traditional approaches. We refer the reader to \cite{LMPC2,norouzi2023integrating} for a more comprehensive discussion of the matter.

Here, we limit ourselves to the study of two different approaches. The first one based on a supervised learning MPC, where the control is learned in an off-line phase, over training data generated by sampling optimal trajectories of agents \cite{albi2021gradient,albi2022supervised,DKS23}. 
The second strategy is based on a sequential linearization, where the control is obtained via Riccati-type equations, and re-calibrated dynamically based on a predictive horizon triggered by decay estimate of the agent-system moments \cite{ahks22} in a similar spirit of event-based MPC \cite{6315560}.
Both approaches are based on the realisation of an efficient control for the agent-based dynamics, however the methods are inherently different: the first strategy aims to approximate the optimal control of non-linear dynamics with a bottom-up procedure, while the second strategy determines a minimum intervention to steer non-linear systems towards a desired goal.

The chapter is organized as follows. In Section 2 we introduce the general setting for optimal control of collective dynamics, including first-order optimality conditions and dynamic programming approaches. Section 3 discusses the aforementioned optimal control techniques in the framework of supervised learning, as a method that enables real-time feedback control. Section 4 presents an alternative approach, moment-driven predictive control, which is based on sequential linearization of the dynamics triggered by a macroscopic observable of the particle ensemble. Section 5 concludes with a series of numerical tests illustrating the performance of the different techniques for large-scale particle systems.

\section{Optimal control for collective dynamics}\label{sec:oc_collectivedyn}

Agent-based models are a natural modelling framework to describe behaviour of a group of interacting entities obeying simple rules \cite{motsch2014heterophilious,carrillo2021controlling}. In this context, it is of great interest to explore the emergence of self-organizing patterns and collaborative behaviours among agents such as consensus, flocking, mills, etc.
Applications include the modelling of language evolution, clustering in opinion dynamics, and collective animal behaviour, among many others.
\cite{blondel2009krause,cucker2004modeling,huang2003individual,couzin2005effective,cucker2007emergent}.  

We consider second-order systems of $N$ interacting agents with state variables described by $s_i=\left(x_i, v_i\right)\in \mathbb{R}^d\times\mathbb{R}^d$ for $i=1\dots N$, where $x_i$ and $v_i$ represent the position and the velocity of the $i-$th agent, respectively. Given an initial configuration of the ensemble, the evolution of the agents is governed by 
\begin{equation}
	\begin{aligned}
		&\dot{x}_i(t)=v_i(t), \\
		&\dot{v}_i(t)=\frac{1}{N} \sum_{j=1}^N a\left(\left\|x_j(t)-x_i(t)\right\|\right)\left(v_j(t)-v_i(t)\right),\qquad i=1, \ldots, N,
	\end{aligned}
\end{equation}
where $a \in C^1(\R_+)$ is a non-increasing positive function encoding the interaction kernel. In particular, we focus on kernels of Cucker-Smale type
\begin{equation}\label{kernel}
	a(r)=\frac{K}{\left(1+r^2\right)^\beta}, \quad \beta \geq 0,\; K>0\,.
\end{equation}
In what follows, for the sake of simplicity, we set $\beta=K=1$. A remarkable property of the Cuker-Smale model is the emergence of consensus by isotropic averaging, i.e., convergence towards $$\bar{v}=\frac{1}{N} \sum_{j=1}^N v_j.$$

Consensus convergence depends on the interaction strength, encoded in $a(\cdot)$, and the cohesiveness of the initial configuration \cite{HaHaKim}. Whenever consensus is not guaranteed to emerge, an external control action can be introduced to enforce convergence towards the desired pattern. To this end, the controlled collective evolution can be modeled as
\begin{equation}\label{eq:csdyn-control}
	\begin{aligned}
		&\dot{x}_i(t)=v_i(t), \\
		&\dot{v}_i(t)=\frac{1}{N} \sum_{j=1}^N a\left(\left\|x_j(t)-x_i(t)\right\|\right)\left(v_j(t)-v_i(t)\right)+u_i(t),
	\end{aligned}
\end{equation}
where each agent is influenced by an ad-hoc control component, with $\mathbf{u}(t):=(u_1,\dots,u_N)^\top$ representing a dynamical external force applied in a centralized manner. This external influence can be tuned by solving an optimal control problem within a set of admissible control signals $\mathcal{U}$: 
\begin{equation}\label{cscost}
	\mathbf{u}^*(t) = \underset{\mathbf{u}(\cdot) \in \mathcal{U}}{\text{argmin}} \,\bigg\{\mathcal{J}(\mathbf{u}(\cdot)):=\int\limits_0^{T} \frac{1}{N}\sum_{j=1}^N\left(\left\|\bar{v}(t)-v_j(t)\right\|^2+\gamma\left\|u_j(t)\right\|^2 \right)d t\bigg\},
\end{equation}
subject to the dynamics \eqref{eq:csdyn-control}. In such a way, the control action is designed to penalize both the distance of the agent's velocities from the target $\Bar{v}$ and the control energy required to steer the system towards the consensus manifold, with parameter $\gamma>0$. 

Necessary optimality conditions for the optimal control problem \eqref{eq:csdyn-control}-\eqref{cscost} are derived from Pontryagin's Maximum Principle (PMP) \cite{Bailo_2018} leading to a two-point boundary value problem for the state variables $s_i$, the associated adjoint variables $r_i=(p_i,q_i)$, and the control signals $u_i$:
\begin{equation}\label{PMP_CS}
	\begin{cases}
		\Dot{x}_i &= v_i\\
		\Dot{v}_i &= \frac{1}{N} \sum\limits_{j=1}^N \dfrac{v_i-v_j}{1 + \|x_i-x_j\|^2} + u_i\\
		\Dot{p}_i &= -\frac{1}{N} \sum\limits_{j=1}^N \dfrac{2(x_i-x_j)}{(1+\|x_j-x_i\|^2)^2}\langle q_j - q_i, v_j-v_i \rangle  \\
		\Dot{q}_i &=-\frac{1}{N}\sum\limits_{j=1}^N \dfrac{q_j-q_i}{1 + \|x_i-x_j\|^2} - \frac2N (\Bar{v}-v_i) + p_i
	\end{cases}
\end{equation}
together with boundary conditions given by $(x_i(0),v_i(0)) = s^0_i$ and $(p_i(T),q_i(T)) = (0,0)$ for all $i=1,\dots,N$.  The PMP system is closed with the optimality condition
\begin{equation}\label{opt_cond_CS}
    \mathbf{u}^*(t) = \underset{u}{\operatorname{argmin}} \bigg(\sum\limits_{i=1}^N \langle q_i,\Dot{v}_i \rangle + \frac{\gamma}{N}\|u_i\|^2_2\bigg) = -\frac{N}{2\gamma}\mathbf{q}^\top\,,
    \end{equation}
    where $\mathbf{q} = (q_1,\dots q_N)^\top$.
It is possible to approximate the optimality conditions above by following a reduced gradient approach \cite{optPDE_Kunisch,Bailo_2018,azmi2020optimal}. Starting with an initial guess for the control signal over the whole horizon, this method iterates forward-backward passes over \eqref{PMP_CS}, while updating the control signal $u(\cdot)$ using \eqref{opt_cond_CS} as a gradient step. Unfortunately, such a method becomes unsuitable for real-time control as the number of agent increases.

\subsection{Dynamic programming and the HJB PDE}
The problem at hand can be written in a more general form as a finite-horizon optimal control problem for a system of $N$ agents:
\begin{equation}\label{ocp1}
	\min _{\mathbf{u}(\cdot) \in \mathcal{U}} \mathcal{J}(\mathbf{u}(\cdot);t_0,\mathbf{s}):=\int\limits_{t_0}^{T} \ell( \mathbf{y}(t))+\frac{\gamma}{N}\|\mathbf{u}(t)\|_2^2 d t,\qquad \gamma>0,
\end{equation}
subject to nonlinear, control-affine constraint
\begin{equation}\label{dyn}
	\dot{\mathbf{y}}(t)=f(\mathbf{y}(t))+g(\mathbf{y}(t)) \mathbf{u}(t), \qquad \mathbf{y}(t_0)=\mathbf{s}\,,
\end{equation}
where the initial condition of the particle ensemble is denoted $\mathbf{s}=(s_1, \ldots, s_N)^{\top}$,  with $s_i \in \R^{2d}$, hence the total dimension of the particle system is  $n=2dN$. We define the control variable as $\mathbf{u}(\cdot) \in \mathcal{U}:=\{\mathbf{u}(t): \mathbb{R}_{+} \rightarrow \mathbb{R}^m$ is measurable$\}$. We further assume the running cost $\ell:\mathbb{R}^{n}\rightarrow\mathbb{R}$ is continuously differentiable and  consider dynamics $f: \mathbb{R}^{n} \rightarrow \mathbb{R}^{n}$, $g: \mathbb{R}^{n} \rightarrow \mathbb{R}^{n \times m}$ in  $\mathcal{C}^1(\mathbb{R}^{n})$. 

We solve the optimal control problem (\ref{ocp1}-\ref{dyn}) by means of dynamic programming. For this, we define the value function of the control problem
\begin{equation}\label{V1}
	V(t,\mathbf{s}):=\inf _{\mathbf{u}(\cdot) \in \mathcal{U}} \mathcal{J}(\mathbf{u}(\cdot);t,\mathbf{s}),\qquad t\in[0,T]\,,
\end{equation}
which, in the finite-horizon case, satisfies the following Hamilton-Jacobi-Bellman PDE
\begin{equation}\label{HJB1}
\left\{\begin{array}{l}\partial_t V(t, \mathbf{s})-\frac{N}{2 \gamma} \nabla V(t, \mathbf{s})^{\top} g(\mathbf{s}) g^{\top}(\mathbf{s}) \nabla V(t, \mathbf{s})+\nabla V(t, \mathbf{s})^{\top} f(\mathbf{s})+\ell(\mathbf{s})=0 \\ V(T, \mathbf{s})=0\end{array}\right.
\end{equation}
After solving the HJB PDE above, it is possible to retrieve the optimal control in terms of the gradient of $V(\cdot)$ :
\begin{equation}\label{opt_cond_HJB1}
\mathbf{u}^*(t, \mathbf{s})=\underset{\mathbf{u} \in \mathbb{R}^m}{\operatorname{argmin}}\left\{\frac{\gamma}{N}\|\mathbf{u}\|_2^2+\nabla V(t, \mathbf{s})^{\top}g(\mathbf{s}) \mathbf{u}\right\}=-\frac{N}{2 \gamma} g^{\top}(\mathbf{s}) \nabla V(t, \mathbf{s}),
\end{equation}
In this way, the control design is in feedback form $\mathbf{u}^* = \mathbf{u}^*(t,\mathbf{s}(t))$, associating an optimal action to every state.
However, equation \eqref{HJB1} is a fully nonlinear, first-order PDE lacking a generally applicable explicit solution. Moreover, the HJB equation is cast over the state space of the system dynamics, which can potentially be of arbitrarily high dimensions -- especially in the context of multi-agent systems, in which the dimensionality scales with both the number $N$ of agents and the dimension $d$ of the underlying state space. Following previously introduced notation, the ensemble state $\mathbf{s}$ encompasses $2dN$ entries.  Bellman himself coined the term \textit{curse of dimensionality} to describe the challenges associated to the computational complexity of solving dynamic programming equations, which becomes particularly evident in large-scale interacting particle systems. In this context, we will discuss computational approaches to circumvent the solution of high-dimensional HJB equations.

\subsection{The link between PMP and HJB}\label{subsec:PMP}
In this section we investigate the connection between Pontryagin Maximum Principle and Dynamic Programming for optimal control design. The approaches are inherently different, as PMP yields open-loop controls satisfying a necessary first order optimality system, whilst solutions of the HJB equation are associated with feedback laws that meet both necessary and sufficient optimality conditions. 
The finite-horizon optimal control problem (\ref{ocp1}-\ref{dyn}) is associated with Hamiltonian 
\begin{equation}
	H(t, \mathbf{s}, \mathbf{u}, \mathbf{r})=\bigg(f(\mathbf{s}(t))+g(\mathbf{s}(t)) \mathbf{u}(t)\bigg)^{\top}\mathbf{r}+\ell\big(\mathbf{s}(t)\big)+\frac{\gamma}{N}\|\mathbf{u}(t)\|_2^2 \,,
\end{equation}
which is defined as a function of ensemble state $\mathbf{s}$, control $\mathbf{u}$, and adjoint variable $\mathbf{r}:=(\mathbf{p},\mathbf{q})$, appearing in this formulation as a time-dependent Lagrange multiplier.  

Pontryagin's principle states requirements that must hold along any optimal trajectory $(\mathbf{s}^*,\mathbf{u}^*,\mathbf{r}^*)$, formulated as the following two-point boundary value problem (TPBVP)
\begin{equation}\label{PMP}
	\begin{cases}
		\dot{\mathbf{s}}^*(t) &= \partial_\mathbf{r} H(t, \mathbf{s}^*, \mathbf{u}^*, \mathbf{r}^*) \\
		\dot{\mathbf{r}}^*(t) &= -\partial_\mathbf{s} H(t, \mathbf{s}^*, \mathbf{u}^*, \mathbf{r}^*) \\
		\mathbf{s}^*(0) &= \;\mathbf{s}_0\,, \quad \mathbf{r}^*(T) \;=\;0\,,
	\end{cases}
\end{equation}
and closed with the first order optimality condition for the optimal control:
\begin{equation}\label{PMP_close}
	\mathbf{u}^*(t) =  \partial_\mathbf{u} H(t, \mathbf{s}^*, \mathbf{u}^*, \mathbf{r}^*) = \frac{N}{2\gamma} g^{\top}\big(\mathbf{s}^*(t)\big) \mathbf{r}^*(t)\,.
\end{equation}
Note that this procedure only provides necessary conditions for optimality, meaning that PMP solutions identify stationary points of $\partial_\mathbf{u} H$, 
with only some of them being global minimizers.  Triples $(\mathbf{s}^*, \mathbf{u}^*, \mathbf{r}^*)$ solving (\ref{PMP}-\ref{PMP_close}) are called \emph{extremal trajectories}, and are computed from  a fixed initial configuration $\mathbf{s}_0$, hence they are not in feedback form.

The connection between Pontryagin's and Bellman's principles in finite horizon settings has been investigated since the early works of analysis in optimal control \cite{pontryagin1962}.
Under suitable assumptions on the regularity of the HJB solution (here $V\in\mathcal{C}^1$), it can be shown that the forward-backward dynamics derived in PMP serve as characteristic curves for the HJB equation \cite{zhou1990maximum,subbotina2006method}. 
Consequently, the value function associated with the optimal control problem can be determined at a specific state-space point $\mathbf{s}$ by solving the TPBVP \eqref{PMP} from the initial condition $\mathbf{s}_0$ and integrating along the resulting optimal trajectory, its gradient being the adjoint variable \cite{barron1986,clarke_vinter,bardicapuzzodolcetta}:
\begin{equation}
    V(t,\mathbf{s}_0) = \int\limits_{t_0}^T \ell(\mathbf{s}^*)+\frac{\gamma}{N}\|\mathbf{u}^*(t)\|_2^2 dt\,,\qquad \nabla V(t,\mathbf{s}_0) = \mathbf{r}^*(t)\,.
\end{equation}

A typical objective in MPC is asymptotic stabilization, which is expressed as an \emph{infinite} horizon optimal control problem, 
to which the previous interpretation does not readily extend as the final condition for the adjoint variable is pathologically set at infinity. Several methods have been proposed in the literature to let $T\to\infty$ within PMP \cite{Frankowska2018}; in this present work, we consider a finite horizon $T$ \emph{long enough}  \cite{Lars1} for the system to reach consensus. In particular, we rely on long-horizon open-loop solves departing from the current system configuration $\mathbf{s}$, and retrieve a control action of feedback nature by associating the first time instance of the control signal to the state $\mathbf{s}$ \cite{azmi2020optima}. As investigated in \cite{relax_DP}, this receding horizon strategy can be used as a relaxation of the Dynamic Programming approach. 

\subsection{State Dependent Riccati Equation}
An alternative approach for solving the infinite horizon optimal control problem is to approximate the solution of the associated dynamic programming HJB equation. We consider a quadratic cost
\begin{equation}\label{ocp}
    \min _{\mathbf{u}(\cdot) \in \mathcal{U}} \mathcal{J}(\mathbf{u}(\cdot), \mathbf{y}):=\int\limits_0^{\infty} \mathbf{s}^{\top}(t) Q \mathbf{s}(t)+\mathbf{u}^{\top}(t) R \mathbf{u}(t) d t,
\end{equation}
subject to control-affine dynamical constraints (\ref{dyn}). The matrix $
Q \in\mathbb{R}^{n \times n}$ and $R\in\mathbb{R}^{m \times m}$ are assumed to be symmetric, such that $Q\succeq 0$,  $R\succ 0$. We further assume $f(0)=0, g(0)=0$. 
Accordingly, the optimal value function 
\begin{equation}\label{V}
	V(\mathbf{s}):=\inf _{\mathbf{u}(\cdot) \in \mathcal{U}} \mathcal{J}(\mathbf{u}(\cdot), \mathbf{s})
\end{equation}
satisfies the following Hamilton-Jacobi-Bellman equation
\begin{equation}\label{HJB}
	\nabla V(\mathbf{s})^{\top} f(\mathbf{s})-\frac{1}{4} \nabla V(\mathbf{s})^{\top} g(\mathbf{s}) R^{-1} g(\mathbf{s})^{\top} \nabla V(\mathbf{s})+\mathbf{s}^{\top} Q \mathbf{s}=0
\end{equation}
with optimal feedback control 
\begin{equation}\label{opt_cond_HJB}
	\mathbf{u}(\mathbf{s})=-\frac{1}{2} R^{-1} g(\mathbf{s})^{\top} \nabla V(\mathbf{s})\,.
\end{equation}
To circumvent directly solving the HJB PDE \eqref{HJB} for $V$, we  resort to the State Dependent Riccati Equation (SDRE) approach, providing suboptimal feedback laws that ensure local stabilibity.
The SDRE method can be applied for system dynamics represented in semi-linear form
\begin{equation}\label{semilin_dyn}
	\dot{\mathbf{s}}=A(\mathbf{s}) \mathbf{s}+B(\mathbf{s}) \mathbf{u}(t)\,,
\end{equation}
provided that the pair $(A(\mathbf{s}),B(\mathbf{s}))$ is stabilizable for all $\mathbf{s}\in\mathbb{R}^n$ \cite{AllaKaliseSimoncini}. 
The SDRE approach is based on the idea that the value function \eqref{V} is non-negative by construction, and with no loss of generality it can be represented as a quadratic form $V(\mathbf{s}) = \mathbf{s}^\top \Pi(\mathbf{s}) \mathbf{s}$, with $\Pi(\cdot)\in\mathbb{R}^{n\times n}$ a symmetric matrix-valued function. By approximating $\nabla V(\mathbf{s}) = 2 \Pi(\mathbf{s}) \mathbf{s}$ and  inserting it into the HJB equation \eqref{HJB}, we obtain the State-dependent Riccati Equation
\begin{equation}\label{sdre}
	A^{\top}(\mathbf{s}) \Pi(\mathbf{s})+\Pi(\mathbf{s}) A(\mathbf{s})-\Pi(\mathbf{s}) B(\mathbf{s}) R^{-1} B^{\top}(\mathbf{s}) \Pi(\mathbf{s})+Q=0
\end{equation}
to be solved for $\Pi(\mathbf{s})$. The resulting feedback law then reads
\begin{equation}\label{opt_cond_SDRE}
	\mathbf{u}(\mathbf{s})=-R^{-1} B^{\top}(\mathbf{s}) \Pi(\mathbf{s}) \mathbf{s}\,.
\end{equation}
The SDRE feedback law is suboptimal compared to the solution of the HJB equation \eqref{HJB}, as the semi-linearization $f(\mathbf{s})=A(\mathbf{s})\mathbf{s}$ is not unique. More substantially, equation \eqref{sdre} holds for the approximation of $\nabla V(\mathbf{s})\approx 2\Pi(\mathbf{s})\mathbf{s}$, whereas a chain rule calculation shows 
\begin{equation}
    \big[\nabla V(\mathbf{s})\big]_k = 2\Pi(\mathbf{s})\mathbf{s} + \sum\limits_{i,j}s_is_j\partial_{s_k}\Pi_{i,j}(\mathbf{s})\,.
\end{equation}
 These two aspects are in no way independent, as one could optimize the choice of $A(\cdot)$ in such a way that minimizes the misfit between the SDRE and the HJB feedback laws \cite{AstolfiSDRE,LucaSDRE}. 

\subsubsection{Frozen Riccati Approach}
The state dependency of the operator $\Pi(\mathbf{s})$ in \eqref{sdre} suggests that the computational complexity associated with the SDRE solution is indeed similar to the original HJB \eqref{HJB}. 
We follow a similar approach as in \cite{BanksFrozenRiccati}, realizing the feedback design in a model predictive control fashion: given the current configuration $\mathbf{s}$ of the system, we freeze the operator $\Pi(\mathbf{s}) \equiv\Pi\in\R^{n\times n}$, meaning that \eqref{sdre} reduces to its algebraic form, where all the state dependencies are neglected by accordingly evaluating all the operators at $\mathbf{s}$
\begin{equation}\label{are}
	\Pi = Q + A^{\top}\Pi\,A - A^{\top}\Pi\,B(R + B^{\top}\Pi\,B)^{-1}B^{\top}\Pi\,A\,.
\end{equation}
The system is then evolved for a short time, after which the current state $\mathbf{s}$ is updated. 
This procedure requires subsequently solving the linear quadratic problem associated with the frozen system, which is addressed via the linear quadratic regulator routine.

\subsubsection{Semilinearization of the Cucker-Smale consensus problem} The semi-linear formulation of the dynamics of interest is derived for the ensemble configuration $\mathbf{s}$, where a two-dimensional agent's state is denoted as $(x,v)=(x^1,x^2,v^1,v^2)$:
\begin{equation}\label{state_transform}
	\mathbf{s} : = \bigg(x^1_1,\,x^2_1\,,\dots\,,x^1_N,\,x^2_N,\,v^1_1,\,v^2_1\,,\dots\,,v^1_N,\,v^2_N\bigg)^{\top}\,,
\end{equation}
and it holds as in \eqref{semilin_dyn} for operators: $B = (\mathbb{O}_{2n},\mathbb{I}_{2n})^{\top}$,
\begin{equation}
	A(\mathbf{s}) = \begin{bmatrix}
		\mathbb{O}_{2n} & \mathbb{I}_{2n}\\
		\mathbb{O}_{2n} &  \Tilde{A}(\mathbf{s})
	\end{bmatrix}\qquad \Tilde{A}_{i,j}(\mathbf{s}) = \begin{cases}
		a(\|x_j-x_i\|) & \text{if $i$ odd, } j = i+1\\
		a(\|x_j-x_i\|) & \text{if $i$ even, } j = i-1\\
		- a(\|x_j-x_i\|) & \text{if } i=j\\
		0 & \text{otherwise}
	\end{cases}
\end{equation}
where $\mathbb{O}_{2n},\mathbb{I}_{2n}$ indicate the zero and the identity matrices in $\R^{2n\times2n}$ respectively, $n=dN$, and $a$ Cucker-Smale type kernel defined in \eqref{kernel}.
Furthermore, aiming at writing the consensus cost functional as in \eqref{ocp}, we notice that 
\begin{equation}
	\sum\limits_{i=1}^{N} \|v_i - \Bar{v}\|^2 = \|\mathbf{v}-C\mathbf{v}\|^2 
\end{equation}
for $\mathbf{v} = (v_1,\dots v_N)^\top$ and a matrix $C\in\R^{n\times n}$ full of zeros, except for the $N\times N$ diagonal blocks, which are valued $N^{-1}\cdot\mathbbm{1}_{N}$, for $\mathbbm{1}_{N}$ denoting a matrix in $\R^{N\times N}$ full of ones.
It follows that
\begin{equation}
	\|\mathbf{v}-C\mathbf{v}\|^2 = \mathbf{v}^{\top}\bigg(\mathbb{I}_{n} + C^{\top}C - 2C\bigg)\mathbf{v}\,.
\end{equation}
Thus, we can write the cost \eqref{cscost} in the form \eqref{ocp} w.r.t. linear operators 
\begin{equation}\label{cost_operators}
	Q = \mathbb{I}_{n} + C^{\top}C - 2C\,,\qquad R = \dfrac{\gamma}{N}\cdot \mathbb{I}_{n}\,.
\end{equation}

\section{Supervised learning MPC approach}\label{SLMPC}
When working with high-dimensional problems, both of the PMP and SDRE approaches for control synthesis may still require considerable CPU times to be implemented in the receding horizon framework. To alleviate this issue, we rely on data-driven supervised learning methods to provide an offline approximation of the feedback law, enabling real-time control.

In the context of regression analysis, our primary goal is to approximate the control action by instructing the model through dataset training, effectively establishing a mapping from system configurations to the target variable of interest, which is either the value function or the feedback control.

Our discussion will  proceed as follows: we will start by addressing the generation of synthetic data and the selection of the target for approximation, then introduce the neural network architecture and training, to conclude by exploring its real-time application in a receding horizon fashion.

\subsection{Gradient-augmented regression and data generation}

We aim at providing a model $\tilde{\mathbf{u}}\approx \mathbf{u}$ acting as a surrogate of the feedback law to accelerate the receding horizon control of the system at hand. Being the synthesis of the control action the objective of the learning operation, the first candidate for approximation is the state-to-control map $\mathbf{s}\mapsto \mathbf{u}(\mathbf{s})$.

As an alternative, another quantity of interest is the value function $V$, from where a control can be derived from $\mathbf{u}\approx\nabla V$, according to the optimality conditions in equation \eqref{opt_cond_HJB} and \eqref{PMP_close}. Besides its theoretical significance, a key advantage of selecting $V$ as a learning objective is its property of being a scalar function, hence relatively easy to approximate. 
Note that information about $\nabla V$ is generated for free as a byproduct of both the SDRE and PMP routines. In the SDRE case, this resides in the ansatz $\nabla V (\mathbf{s}) \approx 2 \Pi(\mathbf{s}) \mathbf{s}$, expressed in terms of the solution operator $\Pi(\cdot)$; in PMP instead, the gradient of the value function coincides with the adjoint variable $\mathbf{r}$. Both the learned (and $V$-learned) control designs will be compared learning objectives, depending on the sampling of a suitable training dataset. 

We sample state configurations $\{\mathbf{s}^{(i)}\}_{i=1}^{N_s}$ and couple them with the associated labels $\mathbf{u}^{(i)} := \mathbf{u}(\mathbf{s}^{(i)}),\;V^{(i)} := V(\mathbf{s}^{(i)}),\;\nabla V^{(i)} := \nabla V(\mathbf{s}^{(i)})$:

\begin{itemize}
	\item
 In the PMP framework, we solve the optimality system (\ref{PMP_CS}-\ref{opt_cond_CS}) fixing each sample as initial configuration $\mathbf{y}=\mathbf{s}^{(i)}$. As done in \cite{azmi2020optimal}, once the extremal trajectory $(\mathbf{s}^*,\mathbf{u}^*,\mathbf{r}^*)$ is computed for time $t\in[0,T]$, the variables of interest can be obtained as follows: the feedback control associated with the initial condition $\mathbf{s}^{(i)}$  is $\mathbf{u}^{(i)} = \mathbf{u}^*(t=0)$, the optimal cost-to-go from $\mathbf{s}^{(i)}$ can be retrieved as the cost integrated along the extremal trajectory, and coincides with the value function $V^{(i)} = V(\mathbf{s}^{(i)},t=0)$, 
 its gradient being $\nabla V^{(i)}=\mathbf{r}^*(t=0)$.
	\item 
For the SDRE approach, we solve the Algebraic Riccati Equation  \eqref{are} associated with freezing \eqref{sdre} at each current state $\mathbf{s}^{(i)}$. The target variables can then be retrieved by using the quadratic ansatz for $V^{(i)} = \mathbf{s}^{(i)^\top} \Pi(\mathbf{s}^{(i)}) \mathbf{s}^{(i)}$, $\nabla V^{(i)} = 2 \mathbf{s}^{(i)^\top} \Pi(\mathbf{s}^{(i)})$ and the first order optimality condition for $\mathbf{u}$ in \eqref{opt_cond_SDRE} \cite{albi2021gradient,albi2022supervised}.
\end{itemize}

\subsection{Feedforward Neural Networks}
We limit our the discussion to models belongings to the class of fully connected Feedforward Neural Networks.  These models approximate generic functions $\phi$ by chains of composition of layer maps $l_1, \ldots, l_M$, i.e., $\tilde{\phi} = \phi_\theta \approx \phi$ with
$$
\phi_\theta(z)=l_M \circ \ldots \circ l_2 \circ l_1(z),\qquad z\in\mathbb{R}^n,
$$
during which the information flows from the input nodes to the output ones in a unidirectional path, avoiding any cycles or loops. For $m=1,\cdots,M$, the layer $l_m$ is defined as $l_m(z)=\sigma_m\left(A_m z+b_m\right)$, where
the nonlinear activation function $\sigma_m: \mathbb{R}^n\rightarrow\mathbb{R}$ is applied component-wise to a linear combination of the layer input variable, with $A_m\in\mathbb{R}^{n\times n}$ and $b_m\in\mathbb{R}^n$ being weight matrix and bias vector for the $m$-th  layer respectively. 

Considering a data set $\mathcal{T}=\left\{z^{(i)}, \phi^{(i)}:=\phi\left(z^{(i)}\right)\right\}_{i=1}^{N_s}$, the $\theta=\left\{A_m, b_m\right\}_{m=1}^M$ 
are trainable parameters to be optimized according to a loss function $\mathcal{L}$ in a supervised learning fashion
$$
\theta=\arg\min _\theta \sum\limits_{z\in\mathcal{T}}\mathcal{L}\left(\phi(z), \phi_\theta(z)\right)\,.
$$
In our setting, we rely on the mean squared error (MSE) as measure of the approximation error $\mathcal{L}_0$. In the gradient-augmented regression for $V$, we penalize the discrepancy with respect to the derivative information as well:
\begin{equation}
	\begin{aligned}
		\mathcal{L}_0(\phi, \phi_\theta) &= \sum\limits_{i=1}^{N_s}\dfrac{\|\phi^{(i)}-\phi_\theta(z^{(i)})\|^2}{N_s}\,,\\
		\mathcal{L}_1(\phi, \phi_\theta) &= \mathcal{L}_0(\phi, \phi_\theta) + \mu \mathcal{L}_0(\nabla \phi, \nabla \phi_\theta)\,. 
	\end{aligned}
\end{equation}

The number of layers $M$, the number of neurons per layer (width of layer), the activation functions $\sigma_m(\cdot)$ in the hidden layers and the loss weight $\mu$ are hyper-parameters which need to be optimally tuned according to the model's goodness of fit, so that the trained model adapt properly to new, previously unseen data.

\subsection{MPC with the learned feedback control}
Once the models $\tilde{\mathbf{u}}$ are trained over the synthetically generated data, they can be used as numerical representations of the optimal feedback control law. Defining an uniformly spaced time-discretization with a step size $\delta$ on the time interval $[0,T]$ by $t_h = h\delta$, $h=1,\dots,T/\delta$, we denote $\mathbf{s}^h := \mathbf{s}(t_h)$. Departing from the initial configuration of the system $\mathbf{s}_0 = \mathbf{y}$, we follow a receding horizon strategy for the controlled evolution
\begin{equation}
    \mathbf{s}^{h+1} = \mathbf{s}^h + \delta\bigg(f(\mathbf{s}^h)+g(\mathbf{s}^h)\tilde{\mathbf{u}}(\mathbf{s}^h)\bigg)\,.
\end{equation}
This allows to replace the subsequent calls to optimal control solvers by mere model evaluations, leading to a considerable complexity reduction, as shown in terms of CPU speed-up in the numerical test section.   

\section{Moment-driven predictive control}\label{sec:MPC}

We propose an additional approach to avoid the limitations associated with the synthesis of optimal feedback laws for high-dimensional non-linear dynamics through conventional approaches, like the computationally expensive solution of the Hamilton-Jacobi-Bellman PDE. To his end, we present a sub-optimal feedback-type control through the linearization of the interaction kernel, and  solve the resulting linear-quadratic optimal control problem through a Riccati equation.  The proposed methodology yields a control law for the linear model, which is later  embedded into the non-linear dynamics.  The  proposed design only requires recurrent measurements of the nonlinear state. No continuous estimation of the nonlinear state or the synthesis of a nonlinear feedback is necessary. The approach relies on the fact, that the linearized system could be controlled efficiently using a Riccati formulation.

However, using the linear optimal control within the nonlinear model does not necessarily yield a stabilizing control law, since over time the nonlinear dynamics may be far from the linearization point. Because of the latter, we aim to quantify the impact of this control. This will determine  the number of linearization updates needed to stabilize the nonlinear system. 
The performance is quantified  by estimating the decay of macroscopic quantities such as the second moment of the particle ensemble. 

This framework, called Moment-driven Predictive Control (MdPC) is  based on the work \cite{ahks22}. Using dynamic estimates of the moments decay, a forward error analysis  is performed to estimate the next linearization point.

Moreover, the proposed control strategy is capable of  treating efficiently high-dimensional non-linear control problems due to the Riccati approach employed for the linearized system.

\subsection{Riccati-based open-loop control}\label{subsec:ffa}

We consider a linearization-based approach for the finite horizon minimization problem introduced in Section \ref{sec:oc_collectivedyn}. In this new approach, the non-linear dynamics \eqref{eq:csdyn-control} is linearized for every agent around an  equilibrium $(\tilde x, \tilde v)$. We assume that the communication function $a$ evaluated at the equilibrium is such that $a(0) \equiv \bar p\,,$ for some $\bar p>0$. 
 Introducing the shift $z_i = x_i-\tilde v t, \ w_i = v_i -\tilde v$, we obtain the  linearized system
 \begin{align*}
	\dot{z}_i (t)&= w_i(t),  &z_i (0)= x_i(0), \\
 	\dot{w}_i(t) &= \frac{1}{N}\sum_{j=1}^N \bar{p}\ (w_j(t)-w_i(t))+ u_i(t), &w_i(0)=v_{i}(0)-\tilde v,
 \end{align*}

 For the sake of compactness, and differently from the previous section, here we consider the full states of the interacting agent system as $z(t),w(t),u(t)\in\R^{N\times d}$, i.e. the position and velocity of each $i$-agent are stored by row.
 Hence, we can write in in matrix-vector notation 
 	\begin{equation*}
 		\begin{bmatrix} \dot z (t) \\ \dot w (t)\end{bmatrix}
 		= \begin{bmatrix} \mathbb{O}_N & \mathbb{I}_N \\ \mathbb{O}_N  & P \end{bmatrix} \begin{bmatrix} z \\ w \end{bmatrix}
 		+ \begin{bmatrix} \mathbb{O}_N  \\\mathbb{I}_N \end{bmatrix} u\, ,
 	\end{equation*}
 where the matrix $P\in\R^{N\times N}$ is defined as
 \begin{align} \label{AB}
 	(P)_{ij}=
 	\begin{cases} &\frac{\bar{p}(1-N)}{N},\qquad i=j,\\
 		&\frac{\bar{p}}{N},\qquad\qquad i\neq j,\\
 	\end{cases} \qquad i,j = 1\ldots,N.
 \end{align}
This second-order system  is controllable \cite{herty2018suboptimal} and since it is linear, the optimal control problem is solved using the Riccati equation. We focus on the case where the average velocity in the cost \eqref{cscost} is a fix target, i.e. $\bar w=0$. The exact solution is given as state feedback  by
\begin{align}\label{eq:Riccati_ctrl}
	u(t) = -\frac{N}{\gamma} K_{22}(t)w(t),
\end{align}
where $K_{22}\in\R^{N\times N}$ is associated to the solution of the following differential Riccati matrix-equation
\begin{equation*}
	\begin{aligned}
		-\begin{bmatrix} \dot K_{11} & \dot K_{12} \\ \dot K_{21} & \dot K_{22} \end{bmatrix}
		= &\begin{bmatrix} \mathbb{O}_N  & \mathbb{O}_N  \\ \mathbb{O}_N  & K_{21} + K_{22}P \end{bmatrix}
		+ \begin{bmatrix} \mathbb{O}_N & \mathbb{O}_N  \\ K_{11} + P K_{21} & K_{12} + PK_{22} \end{bmatrix} \\
		&- \frac{N}{\gamma} \begin{bmatrix} K_{12} K_{21} & K_{12} K_{22} \\ K_{22} K_{21} & K_{22}K_{22} \end{bmatrix}
		+ \begin{bmatrix} \mathbb{O}_N  & \mathbb{O}_N  \\ \mathbb{O}_N  & \mathbb{I}_N \end{bmatrix} ,
	\end{aligned}
\end{equation*}
with terminal conditions $K_{ij}(T) = \mathbb{O}_N $, for $\ i,j = 1,2$. This system is easily solved with $K_{11} = K_{12} = K_{21} = \mathbb{O}_N $ and $K_{22}$ fulfilling 
\begin{equation}\label{eq:Riccati}
	- \dot{K}_{22} = K_{22}P+P K_{22}-\frac{N}{\gamma} K_{22}K_{22} + \mathbb{I}_N, \quad K_{22}(T) = \mathbb{O}_N.
\end{equation}
For a general linear system, we need to solve the $N\times N$ differential system \eqref{eq:Riccati}, which can be costly for large-scale agent-based dynamics. However, we exploit the symmetric structure of the matrix $P$ to reduce the dimension of the Riccati equation  considering $K_{22}$ parametrized with a diagonal element $k_d(t)$ and an off-diagonal element $k_o(t)$, according to the approach in \cite{ahks22,herty2015mean}, as follows
\[
(K_{22})_{ij}=\delta_{ij} \, Nk_d(t)+(1-\delta_{ij})N^2\,k_o(t),\qquad i,j=1,\ldots,N,
\]
with $\delta_{ij}$ indicating the Kronecker delta.
Hence, we have the following
\begin{prop} \label{prop:kdko}
	For the linearized dynamics, the solution of the Riccati equation \eqref{eq:Riccati} reduces to the solution of
	\begin{subequations}\label{eq:kd_ko}
	\begin{align}
		-\dot\kd &= -2\bar{p}\alpha(N)\left(\kd - \frac{\ko}{N}\right) - \frac{1}{\nu}\left(\kd^2+\frac{\alpha(N)}{N}\ko^2\right) + 1, \qquad k_d(T)=0, \label{kd}
		\\
		-\dot\ko &=   2\bar{p}\left(\kd - \frac{\ko}{N}\right) - \frac{1}{\nu}\left(2\kd \ko+\alpha(N)\ko^2-\frac{1}{N}\ko^2\right),\qquad k_o(T)=0\,, \label{ko}
	\end{align}	
	\end{subequations}
	with terminal conditions $k_d(T)=k_o(T)=0$, and where $\alpha(N) = \frac{N-1}{N}$.
\end{prop}
The state feedback  \eqref{eq:Riccati_ctrl} is given  by
\begin{align}\label{eq:ctrl_ricc}
	u_i(t) &= - \frac{1}{\gamma}\left(\left(\kd(t)-\frac{\ko(t)}{N}\right) w_i(t) + \frac{\ko(t)}{ N} \sum_{j=1}^N w_j(t)\right)\,.
\end{align}
We refer to \cite{ahks22} for the proof of Proposition \ref{prop:kdko} and further details regarding this approach.

In order to approximate the synthesis of feedback laws for the original non-linear optimal control problem, we study sub-optimal stabilizing strategies induced by the feedback \eqref{eq:ctrl_ricc}. In a feed forward approach we 
apply \eqref{eq:ctrl_ricc} directly in the nonlinear dynamics. 
This  is an open-loop type of control, since all the information on the state of the non-linear system reduces to the initial state of the linearized system, assuming $w_i(0)=v_i(0)$. While it is clearly outperformed by an optimal feedback law in terms of robustness, it has the advantage that it can be implemented without requiring a continuous measurement of the full nonlinear state $v(t)$, making it appealing for systems where recovering the true state of the dynamics can be expensive or time-consuming. This stabilization strategy is clearly suboptimal with respect to the original optimal control problem, and in general will not guarantee the stabilization of the non-linear dynamics \eqref{eq:csdyn-control}.  In order to estimate these performances in the case where a large number of agents is present, i.e. $N\gg 1$, we conduct the error analysis from a {\em mean-field} perspective.

\subsection{Moments estimates for error analysis}\label{sec:mean-field}

We consider the probability density distribution of agents  to describe the collective behavior of a large ensemble of agents, and  retrieve upper and lower bounds for the decay of the mean-field density towards the desired configuration.
For further details on mean-field derivation of particle systems we refer to \cite{ahks22,carrillo2010particle,carrillo2014derivation,canizo2011well} as well as to \cite{MR3264236,MR4028474} for rigorous results on the convergence. We denote by $f(t,z,w)$ and $g(t,x,v)$ the mean-field probability densities corresponding to the linear and non-linear dynamics, respectively. The detailed equations are stated in \cite{ahks22}. The performance of the MdPC approach relies on estimates of the variance of $g$. For given initial data  $g^0:= g(0,x,v)$, we have the following result
\begin{prop}\label{prop:Bounds}
Assume the kernel $a(\cdot)$ to be a bounded function, namely
\begin{align*}
	a(r)\in [-\alpha,\beta],\qquad \alpha,\beta\geq 0.\qquad 
\end{align*}
	We have the following lower and upper bounds for the evolution of the variance $\sigma^2[g]$:
	\begin{equation*}\label{ul_bounds}
		\begin{split}
			\sigma^2[g^0]e^{-2\beta t} \left(1 - B^+_\beta(0,t)\right)^2
			\leq
			\sigma^2[g](t)
			\leq
			\sigma^2[g^0]e^{2\alpha t} \left( 1 + B^-_\alpha(0,t) \right)^2,\quad\text{where}
		\end{split}
	\end{equation*}
	\begin{equation*} \label{B}
		\begin{split}
			B^\pm_c(t_0,t) = \frac{1}{\nu}\int_{t_0}^{t-t_0}  \eta(s-t_0)\kd(s)e^{\pm c (s-t_0)} ds,\cr
			\eta(t-t_0) = \exp\left\{-2 \bar{p}(t-t_0)-\frac{1}{\nu} \int_{t_0}^{t-t_0} \kd(r) dr\right\}.
		\end{split}
	\end{equation*}
\end{prop}

	\subsection{MdPC approach} \label{sec:MdPC}

Relevant issues in the MPC literature are the selection of suitable intermediate time horizons that can ensure asymptotic stability, as well as the design of effective optimization methods, such that the implementation is suitable for real-time control.
We propose a MPC-type algorithm where instead of fixing a prediction horizon, the re-calibration of the control law introduced in Section \ref{subsec:ffa} is triggered adaptively in time. This is based on  a direct estimate of the  decay of the variance estimated above.

Starting by the control \eqref{eq:ctrl_ricc} we consider densities at initial time given by $g(0,x,v)=g^0(x,v)$, $f^0(x,v)\equiv g^0(x,v)$.  
Prediction of the error in the variance decay $\sigma^2[g](t)$ is given by Proposition \ref{prop:Bounds} that yields
\begin{equation*}\label{eq:chk_bound}
	\Delta_{\sigma}(t_0,t) = \sigma^2[g(t_0,x,v)] \left( e^{2\alpha(t-t_0)} \left( 1 + B^-_\alpha(t_0,t) \right)^2 - e^{-2\beta(t-t_0)} (1 - B^+_\beta(t_0,t))^2 \right).
\end{equation*}
Then, $\Delta_{\sigma}(t_0,t)$  is used to guarantee the variance of $g$ is below a fixed threshold $\delta_{tol}>0$.
This defines a time $t_1>t_0$ such that $\Delta_{\sigma}(t_0,t_1)>\delta_{tol}$. Once we reach the threshold, we  reinitialize by updating the state of the linearized dynamics at time $t_1$ by setting $f(t_1,x,v) \equiv g(t_1,x,v)$ and repeat the procedure, see Algorithm \ref{alg:MPCsigma}.
\begin{algorithm}[!ht]
	\caption{[MdPC]}\label{alg:MPCsigma}
	\begin{algorithmic}
		\STATE 0. Set $k \leftarrow 0$, $t_k = 0$, $g^k(x,v) = g(0,x,v)$, $f^k(x,v) = g(0,x,v)$ and tolerance $\delta_{tol}$\;
		\STATE 1. Solve the Riccati equation to obtain $\kd$, $\ko$ on the time interval $[0,T]$\;
		\STATE 2. Find the time $t_{k+1}$ such that $t_{k+1} := \min \lbrace t \vert t_k< t \leq T,\Delta_{\sigma}(t_k,t) > \delta_{tol} \rbrace$
		\WHILE{$t_{k+1}\leq T$}
		\STATE i. Evolve the linear and non-linear dynamics up to $t_{k+1}$\;
		\STATE ii. Set $g_{k+1}(x,v) = g(t_{k+1},x,v)$,  $f_{k+1}(x,v) = g(t_{k+1},x,v)$\;
		\STATE iii. $k \leftarrow k+1$
		\STATE iv. Compute $t_{k+1}$ from step 2\;
		\ENDWHILE
	\end{algorithmic}
\end{algorithm}

\section{Numerical Tests}
We apply the proposed methodologies to the high-dimensional second-order consensus problem \eqref{cscost} for a system of $N = 50$ agents in $\mathbb{R}^2$, following dynamics of Cucker-Smale type \eqref{eq:csdyn-control} and a control penalization parameter $\gamma = 0.1$. We consider the integration of the system in a time-frame $t\in [0,T]$ with $T=10$ and time-step $\delta = 0.01$. The solution of this optimal control problem is associated with an HJB equation cast in a $2dN=200$ dimensional space, hence the need to rely on numerical representations for the optimal feedback design. 

\subsection{Test 1: Supervised Learning MPC}
In this first numerical example, we test the neural network approximation presented in Section  \ref{SLMPC}.
Once the training data are collected, we train models for both the control designs and target variables, denoted as:
\begin{itemize}
	\item $\mathbf{u}_\theta^{SDRE}$, $\mathbf{u}_\theta^{PMP}$ -- for the direct approximation of the mapping $x\mapsto u(x)$;
	\item $\mathbf{u}_V^{SDRE}$, $\mathbf{u}_V^{PMP}$ -- for the approximation of the feedback map through a model for the value function $V$. In this case, a deterministic layer is stacked on top of the FNN, encoding (\ref{PMP_close},\ref{opt_cond_HJB}).
\end{itemize}

We prescribe the initial condition for the ensemble state $\mathbf{s} \in [0,1]^{2dN}$, and a grid search in the parametric space is performed, to select an appropriate model architecture. In Table \ref{tab:specs}, we compare the architectures of the different models, characterized by the number $M$ of layers, the width of each one of them, the activation functions $\sigma$ on the hidden layers and the gradient augmentation parameter $\mu$. All models have been trained with the Adam solver \cite{Adam} over batches of $200$ samples. In the table, we also display the  Coefficient of Variation or Percent Root Mean Squared Error (PRMSE) for the different models when applied to a testing dataset of $10^6$ unseen sample initial conditions in $[0,1]^{2dN}$. 

\begin{table}
	\centering
	\begin{tabular}{cccccc}
		& M & layers width & $\sigma$ & $\mu$ &  PRMSE\\\hline\\[-7pt]
		\(\mathbf{u}_\theta^{SDRE}\)  & $1$ & $1000$ & tanh & $-$ & $0.7921\%$ \\[3pt]
		\(\mathbf{u}_\theta^{PMP}\)  & $3$ & $400$ & sigmoid &$-$& $1.1477\%$ \\[3pt]
		\(\mathbf{u}_V^{SDRE}\)  & $2$ & $200,\,100$ & sigmoid &$0.07$& $4.1843\%$ \\[3pt]
		\(\mathbf{u}_V^{PMP}\)  & $2$ & $800$ & sigmoid &$0.1$& $4.8775\%$\\[3pt]
	\end{tabular}
	\caption{\textbf{Test 1:} Different model specifics, including architectural details and  Percent RMSE in a test-set.}
	\label{tab:specs}
\end{table}

\begin{figure}[t]
	\centering
	\includegraphics[width=0.6\textwidth]{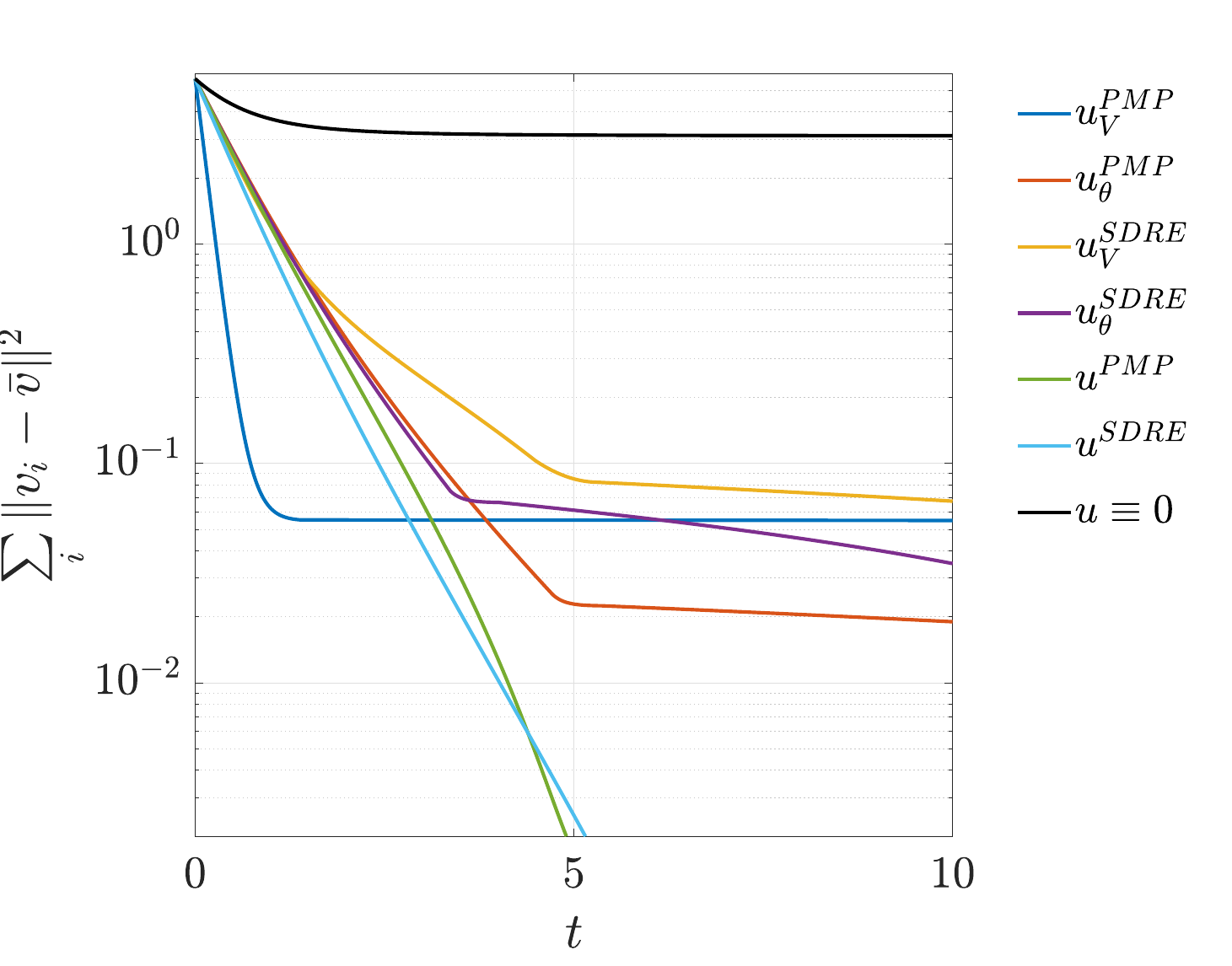}
	\caption{\textbf{Test 1:} Comparison of learned and target trajectories for a random initial condition in $[0,1]^{2dN}$. All the models improve consensus with respect to the uncontrolled dynamics. In terms of learned models, the $\mathbf{u}_V^{PMP}$ model can be considered the optimal one, with an overall cost (at final time horizon $T$) $\mathcal{J}\approx1.44$. }
	\label{fig:dissensus}
\end{figure}

In Figure \ref{fig:dissensus} a plot of the distance from the target configuration is displayed for the different models, and compared with the one resulting from both the PMP approach and the frozen SDRE, in a simulation in times $t\in[0,10]$, with $\delta = 0.01$ seconds.

\begin{figure}[t]
	\centering
	\includegraphics[width=0.329\textwidth]{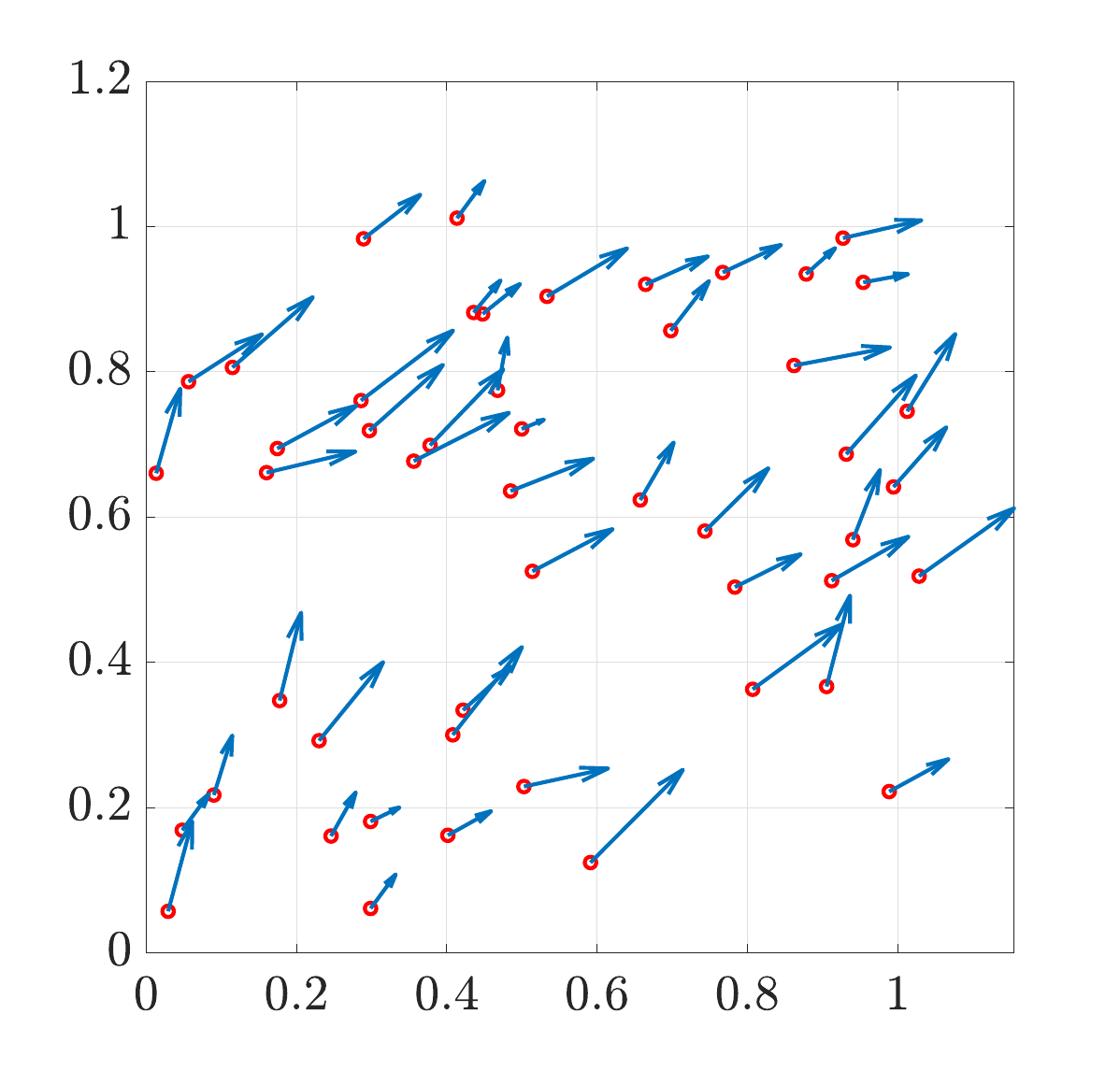}
	\includegraphics[width=0.329\textwidth]{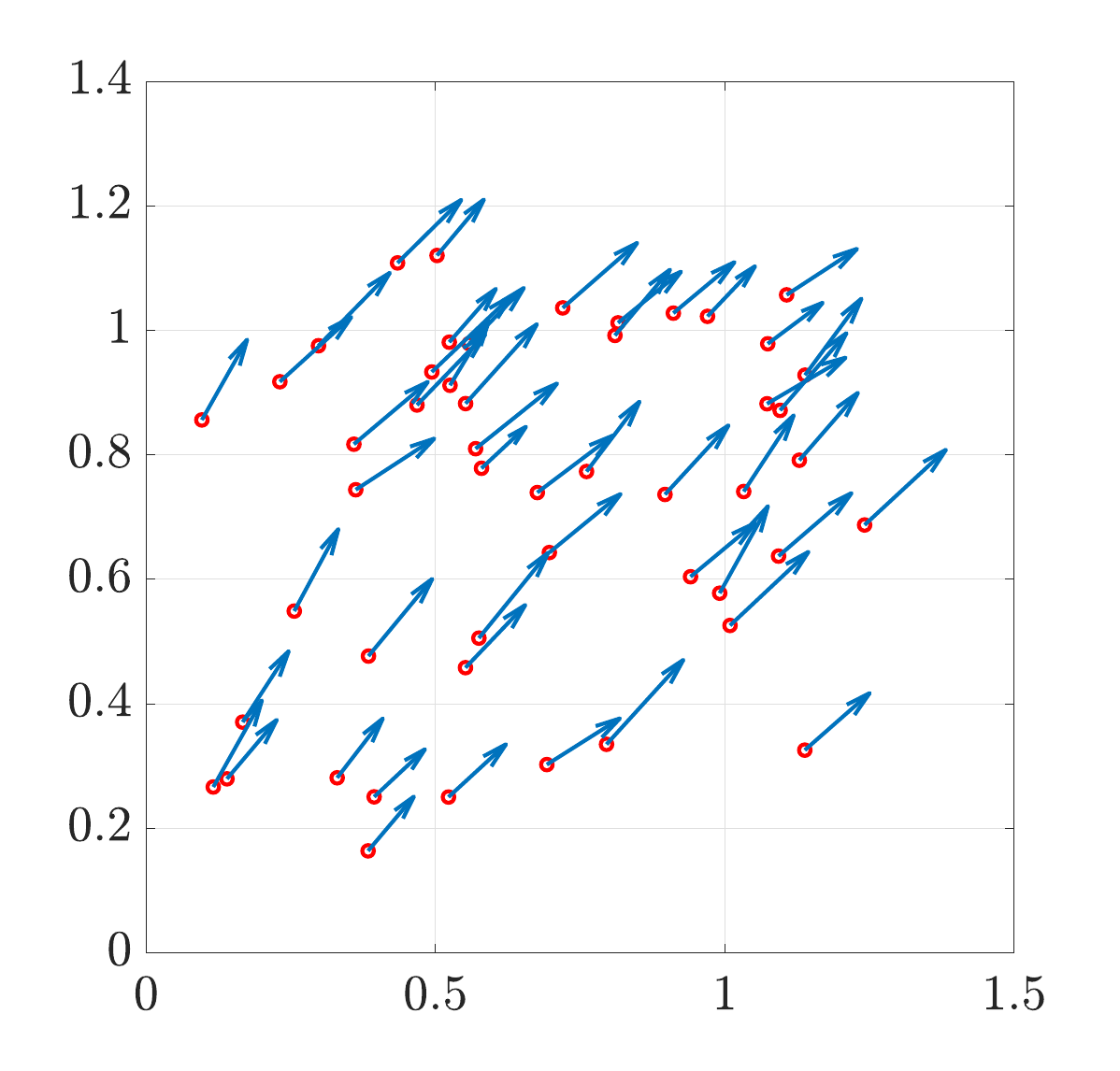}
	\includegraphics[width=0.329\textwidth]{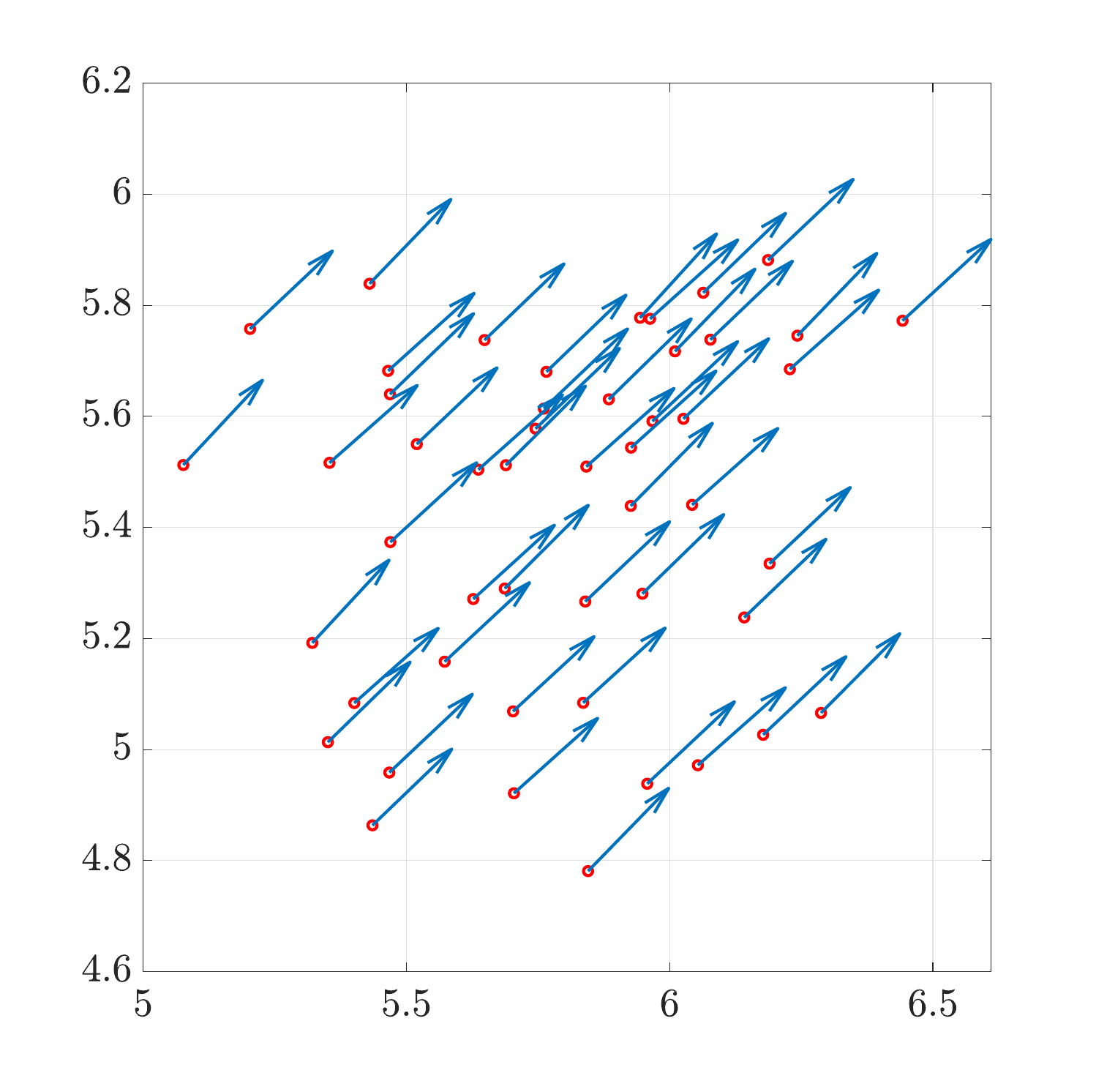}
	\caption{\textbf{Test 1:} System configuration at time $t=0,1,10$ seconds respectively, under the control action provided by the feedback map $\mathbf{u}_V^{PMP}$. }
	\label{fig:microscopic_dyn}
\end{figure}

Another interesting comparison is the one in terms of the required CPU time for integrating the system trajectories while controlling them via the different methods. This is presented in Table \ref{tab:times}, versus the CPU time used by the forward simulation with the original solvers.

\begin{table}
	\centering
	\begin{tabular}{cccccc}
		$\mathbf{u}_\theta^{SDRE}$ & $\mathbf{u}_\theta^{PMP}$ & $\mathbf{u}_V^{SDRE}$ & $\mathbf{u}_V^{PMP}$ & SDRE & PMP\\\hline\\[-7pt]
		$13.9152$ & $9.2503$& $11.8068$  & $4.0699$ & $1821.21$ & $425.1578$ \\
	\end{tabular}
	\caption{\textbf{Test 1:} CPU time (in seconds) required for the system integration along time $t\in[0,10]$ under the different alternative feedback designs.  
	}\label{tab:times}
\end{table}

\subsection{Test 2: MdPC for the Cucker-Smale model}

For the second experiment, we employ the MdPC method outlined in Section \ref{sec:MPC} to address the high-dimensional control problem given by \eqref{cscost} along with the constraint \eqref{eq:csdyn-control}. Using the Cucker-Smale interaction kernel and maintaining the parameters from Test 1, we modify the initial velocities while keeping the positions initialized as before, i.e., $v_i(t)$ uniformly distributed in the interval $[-1,1]$ and $x_i(t)$ in $[0,1]$. The control is determined through the MdPC algorithm \ref{alg:MPCsigma}, and we use Proposition \ref{prop:kdko} to derive the corresponding Riccati equations. We conduct a numerical investigation about the decay of the variance $\sigma^2$ with varying tolerances $\delta_{tol}$.

\begin{figure}[t]
	\centering
	\includegraphics[width=0.329\textwidth]{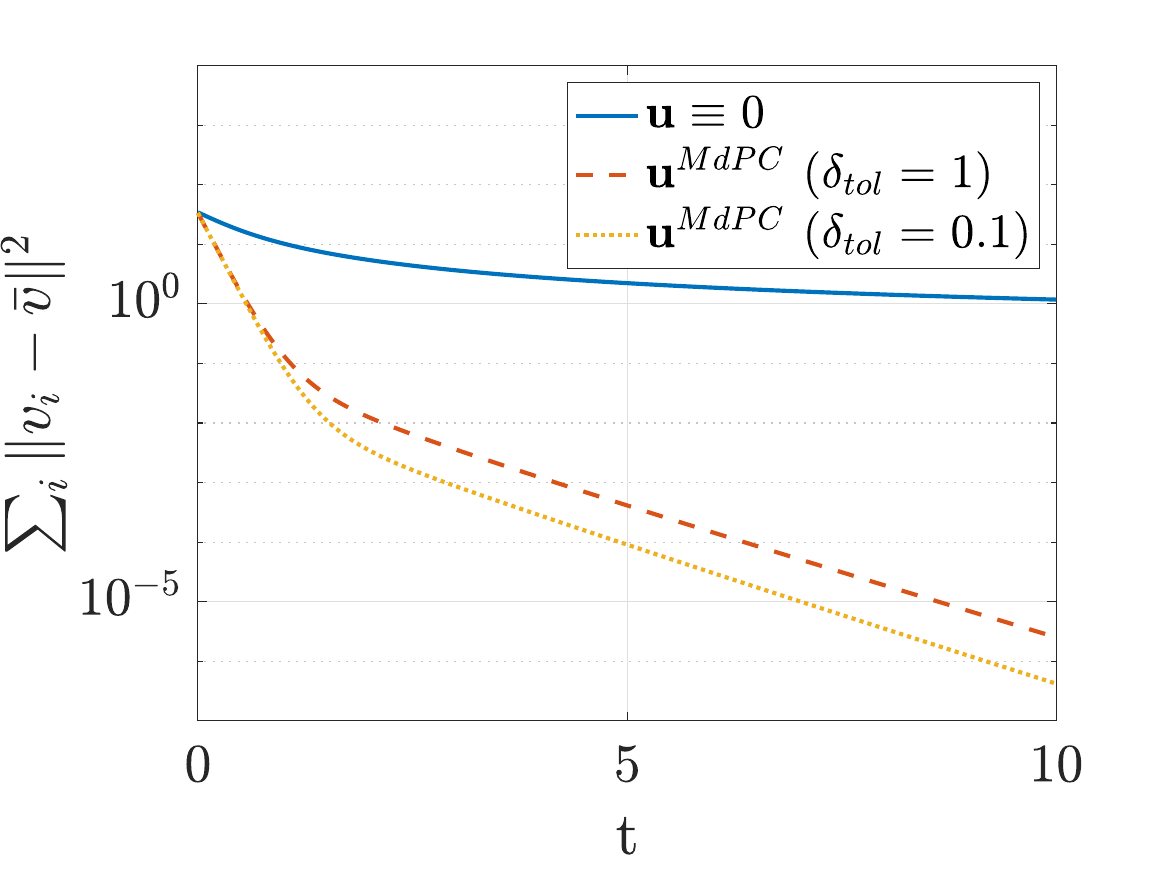}
	\includegraphics[width=0.329\textwidth]{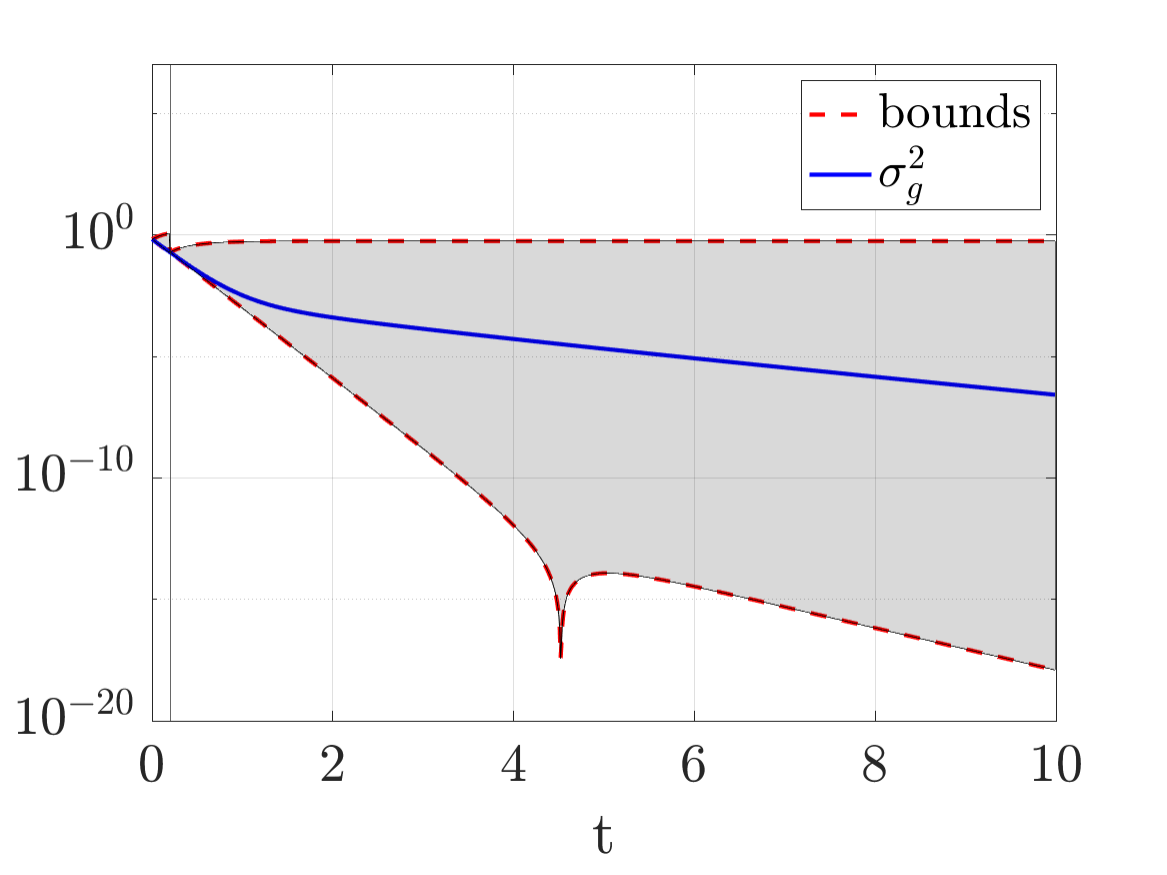}
 	\includegraphics[width=0.329\textwidth]{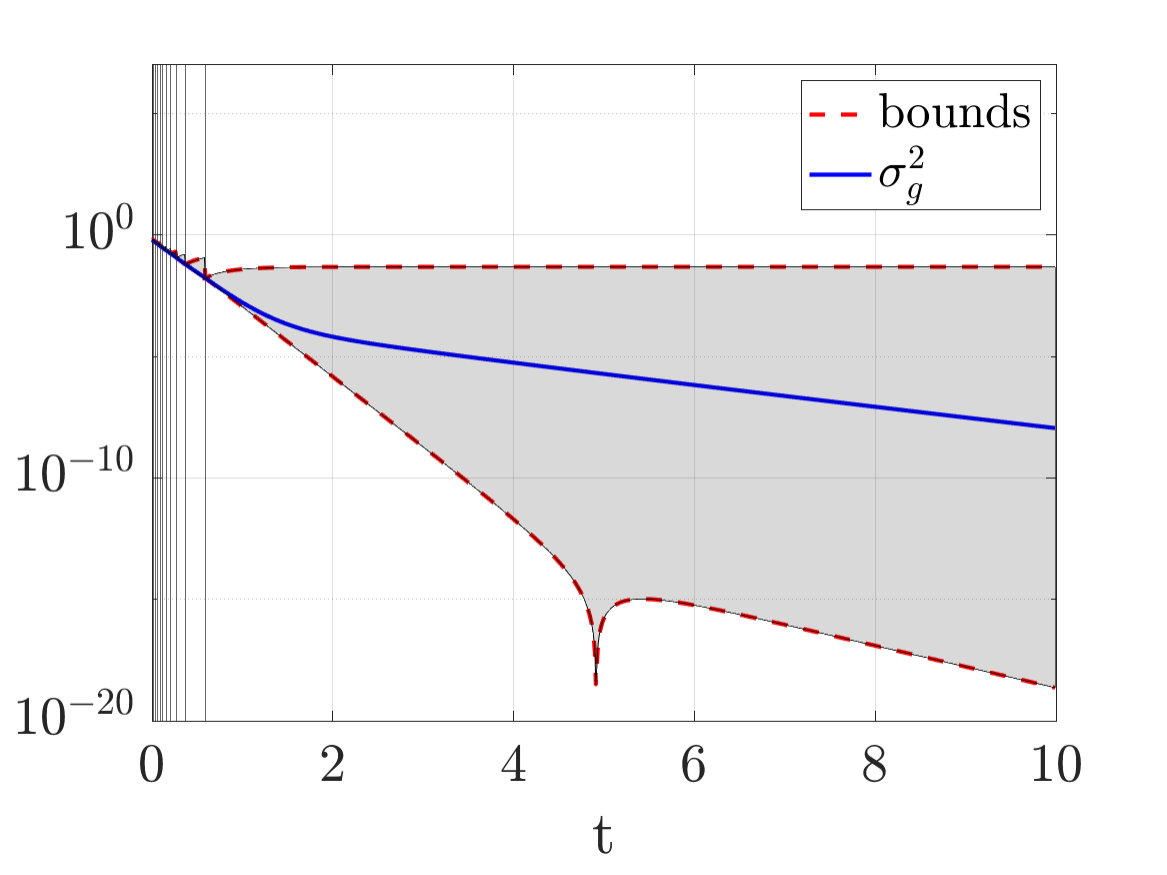}
	\caption{\textbf{Test 2:} On the left, comparison of MdPC-controlled and uncontrolled trajectories in terms of particles' target distance.  Variance decay together with updates and variance bounds in the middle ($\delta_{tol}=1$) and on the right ($\delta_{tol}=0.1$).}
	\label{fig:mdpc_1}
\end{figure}

In the left panel of Figure \ref{fig:mdpc_1}, various control actions are compared based on target distance, as in Figure \ref{fig:dissensus} from Test 1. Notably, the decay is more pronounced in the presence of control and with a smaller tolerance $\delta_{tol}$. This is due to the fact that we have more information about the non-linear dynamics as the tolerance decreases. The center and right panels of Figure \ref{fig:mdpc_1} depict the variance decay for $\delta_{tol} =1$ and $\delta_{tol} =0.1$ respectively. Black vertical lines represent updates, while dashed red lines illustrate the evolution of the variance bounds. The decay is stronger, and the bounds are stricter with a smaller $\delta_{tol}$. For $\delta_{tol} =1$, only one update is required, corresponding to a single-step observation of the non-linear dynamics (refer to Algorithm \ref{alg:MPCsigma}). While with $\delta_{tol} =0.1$, there are 10 updates over a total of 1000 time-steps.

\begin{figure}[t]
	\centering
	\includegraphics[width=0.329\textwidth]{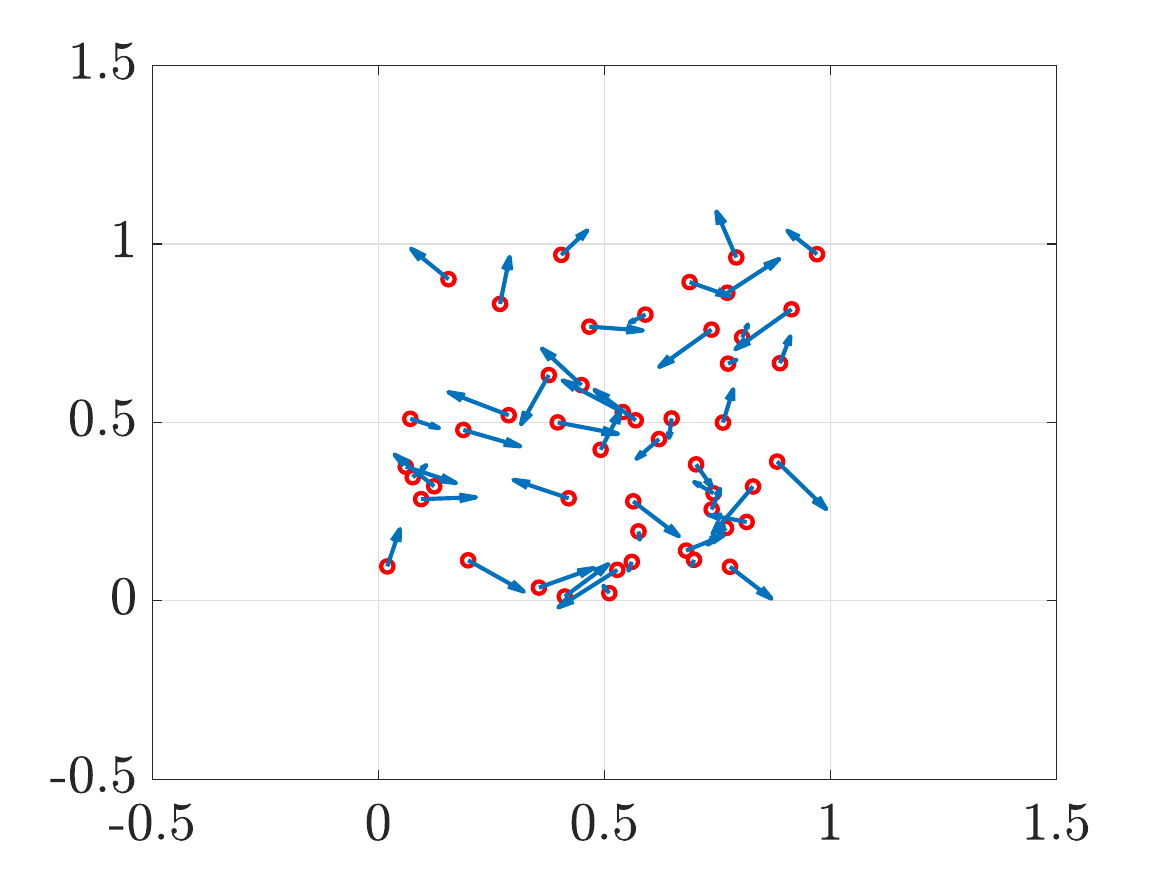}
	\includegraphics[width=0.329\textwidth]{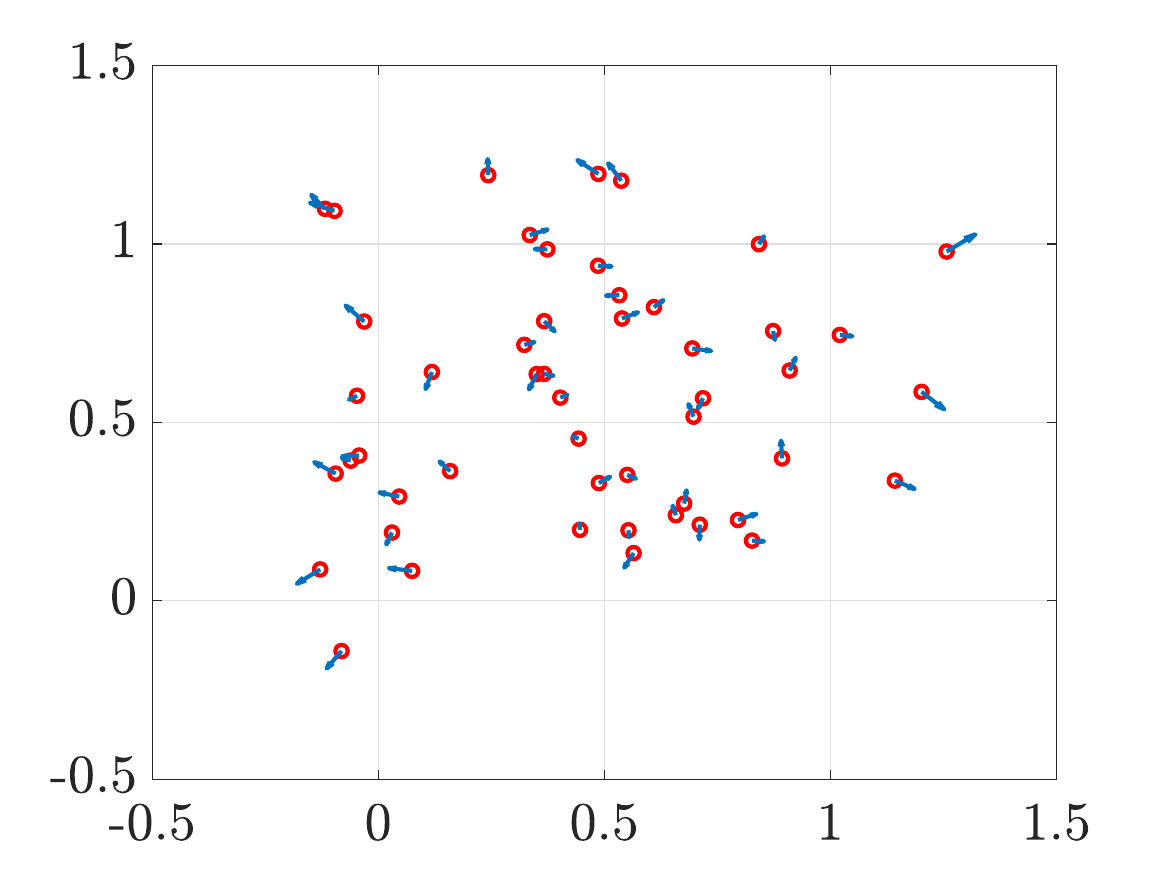}
	\includegraphics[width=0.329\textwidth]{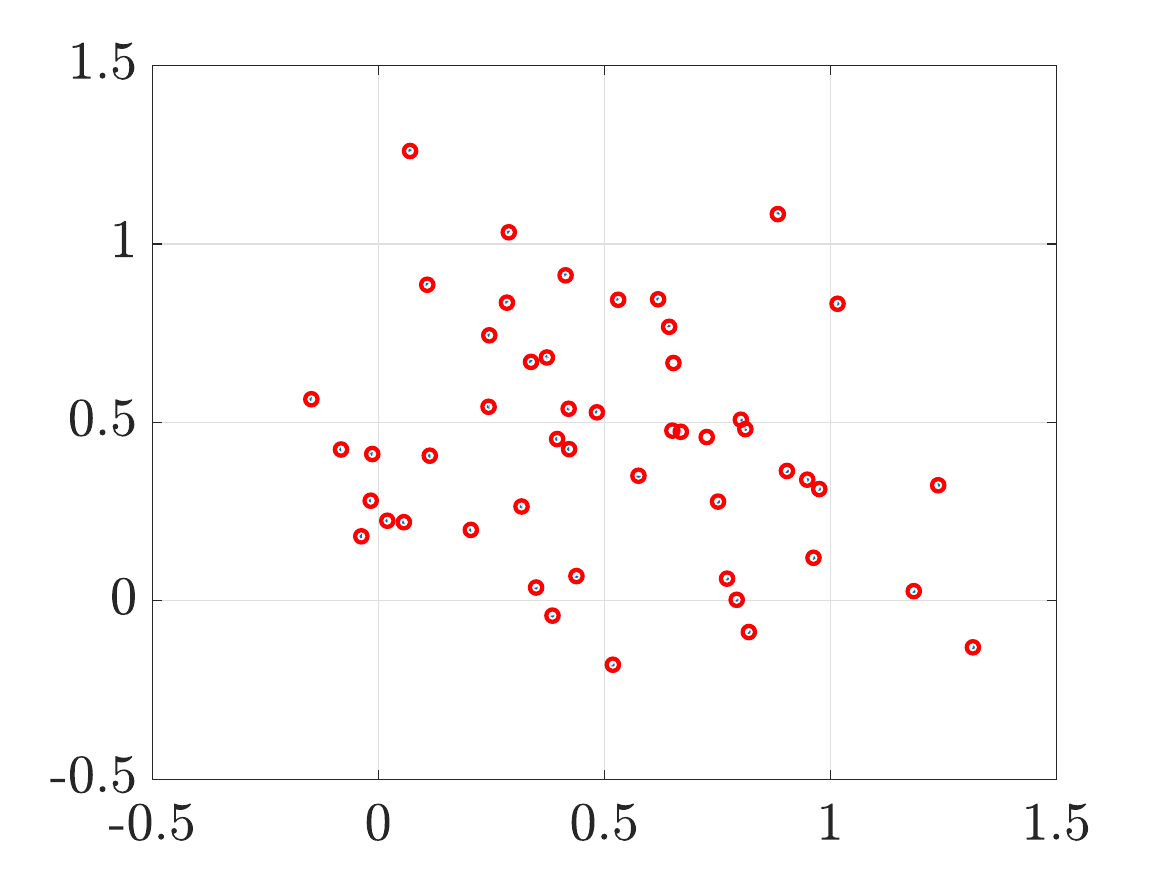}
	\caption{\textbf{Test 2:} System configuration at time $t=0,1,10$ seconds respectively, under the control action provided by $\mathbf{u}^{MdPC}$ with tolerance $\delta_{tol} = 0.1$.}
	\label{fig:mdpc_2}
\end{figure}

In Figure \ref{fig:mdpc_2}, the evolution of positions and velocities is illustrated through three snapshots: the initial data at $t=t_0=0$ on the left, an intermediate time-step at $t=1$ in the center, and the final configuration at $t=T=10$ on the right. The initial mean velocity is zero and the objective is to achieve consensus, hence all final velocities are uniformly zero as expected.

We conclude the analysis by presenting the computational times in Table \ref{tab:MdPC_CPU}.

Using the Riccati equation for control synthesis yields an optimal control for cost functionals associated with linear dynamics only. However, in the more general non-linear case, such control becomes sub-optimal, as evidenced in previous works such as \cite{ACFK17,MR2861587}. We emphasize that the MdPC approach offers an adaptive synthesis of sub-optimal feedback controls specifically designed for non-linear dynamics.
 
\begin{table}
	\centering
	\begin{tabular}{ccc}
		$\mathbf{u}^{MdPC} \ (\delta_{tol} = 1)$ & \hspace{0.5cm} &  $\mathbf{u}^{MdPC} \ (\delta_{tol} = 0.1)$ \\\hline\\[-7pt]
		$0.221$ & &  $0.227$ \\
	\end{tabular}
	\caption{\textbf{Test 2:} CPU time (in seconds) required for the system integration along time $t\in[0,10]$ under the MdPC approach for different tolerances $\delta_{tol}$.
	}\label{tab:MdPC_CPU}
\end{table}

\paragraph{Concluding remarks.} We have reviewed two alternatives to traditional MPC algorithms when controlling large-scale interacting particle systems. Firstly, we have shown that the use of supervised learning techniques in conjunction with deep neural networks is a suitable method for building a suboptimal feedback law based on a finite number of offline samples from optimal trajectories. Synthetic optimal trajectory generation can be done either using first-order optimality conditions in Pontryagin form or by following a State-dependent Riccati Equation approach. In both cases, the high-dimensional control system makes online computation prohibitively expensive, but can be addressed in an offline training phase. Secondly, we propose a moment-driven predictive control framework, where fast online feedback synthesis is achieved by linearization, which is re-computed based on macroscopic quantities of the particle system. This method is suboptimal when compared to nonlinear control laws, but provides a fast computable alternative which do not rely on offline training of a feedback law.

\bibliographystyle{abbrv}
\bibliography{refs,biblioCS}

\begin{thebibliography}{10}

\bibitem{bellomo20review}
G.~Albi, N.~Bellomo, L.~Fermo, S.-Y. Ha, J.~Kim, L.~Pareschi, D.~Poyato, and
  J.~Soler.
\newblock Vehicular traffic, crowds, and swarms: from kinetic theory and
  multiscale methods to applications and research perspectives.
\newblock {\em Math. Models Methods Appl. Sci.}, 29(10):1901--2005, 2019.

\bibitem{albi2021gradient}
G.~Albi, S.~Bicego, and D.~Kalise.
\newblock Gradient-augmented supervised learning of optimal feedback laws using
  state-dependent {R}iccati equations.
\newblock {\em IEEE Control Systems Letters}, 6:836--841, 2021.

\bibitem{albi2022supervised}
G.~Albi, S.~Bicego, and D.~Kalise.
\newblock Supervised learning for kinetic consensus control.
\newblock {\em IFAC-PapersOnLine}, 55(30):103--108, 2022.

\bibitem{MR3542027}
G.~Albi, M.~Bongini, E.~Cristiani, and D.~Kalise.
\newblock Invisible control of self-organizing agents leaving unknown
  environments.
\newblock {\em SIAM J. Appl. Math.}, 76(4):1683--1710, 2016.

\bibitem{albi2023data}
G.~Albi, E.~Calzola, and G.~Dimarco.
\newblock A data-driven kinetic model for opinion dynamics with social network
  contacts.
\newblock {\em prepint arXiv 2307.00906}, 2023.

\bibitem{ACFK17}
G.~Albi, Y.-P. Choi, M.~Fornasier, and D.~Kalise.
\newblock Mean field control hierarchy.
\newblock {\em Appl. Math. Optim.}, 76(1):93--135, 2017.

\bibitem{ahks22}
G.~Albi, M.~Herty, D.~Kalise, and C.~Segala.
\newblock Moment-driven predictive control of mean-field collective dynamics.
\newblock {\em SIAM Journal on Control and Optimization}, 60(2):814--841, 2022.

\bibitem{MR3268062}
G.~Albi, L.~Pareschi, and M.~Zanella.
\newblock Boltzmann-type control of opinion consensus through leaders.
\newblock {\em Philos. Trans. R. Soc. Lond. Ser. A Math. Phys. Eng. Sci.},
  372(2028):20140138, 18, 2014.

\bibitem{ExplicitMPC}
A.~Alessio and A.~Bemporad.
\newblock A survey on explicit model predictive control.
\newblock {\em Nonlinear Model Predictive Control: Towards New Challenging
  Applications}, pages 345--369, 2009.

\bibitem{AllaKaliseSimoncini}
A.~Alla, D.~Kalise, and V.~Simoncini.
\newblock State-dependent riccati equation feedback stabilization for nonlinear
  pdes.
\newblock {\em Advances in Computational Mathematics}, 49:1--32, 2021.

\bibitem{Lars1}
N.~Altmüller and L.~Gr\"une.
\newblock Distributed and boundary model predictive control for the heat
  equation.
\newblock {\em GAMM-Mitteilungen}, 35(2):131--145, 2012.

\bibitem{azmi2020optimal}
B.~Azmi, D.~Kalise, and K.~Kunisch.
\newblock Optimal feedback law recovery by gradient-augmented sparse polynomial
  regression.
\newblock {\em Journal of Machine Learning Research}, 22(48):1--32, 2021.

\bibitem{Bailo_2018}
R.~Bailo, M.~Bongini, J.~A. Carrillo, and D.~Kalise.
\newblock Optimal consensus control of the {C}ucker-{S}male model.
\newblock {\em {IFAC}-{PapersOnLine}}, 51(13):1--6, 2018.

\bibitem{BanksFrozenRiccati}
H.~T. Banks, B.~M. Lewis, and H.~T. Tran.
\newblock Nonlinear feedback controllers and compensators: a state-dependent
  riccati equation approach.
\newblock {\em Computational Optimization and Applications}, 37:177--218, 2007.

\bibitem{bardicapuzzodolcetta}
M.~Bardi, I.~C. Dolcetta, et~al.
\newblock {\em Optimal control and viscosity solutions of
  {H}amilton-{J}acobi-{B}ellman equations, Section 3.4}, volume~12.
\newblock Springer, 1997.

\bibitem{barron1986}
E.~N. Barron and R.~Jensen.
\newblock The pontryagin maximum principle from dynamic programming and
  viscosity solutions to first-order partial differential equations.
\newblock {\em Transactions of the American Mathematical Society},
  298(2):635--641, 1986.

\bibitem{MR2974186}
N.~Bellomo and J.~Soler.
\newblock On the mathematical theory of the dynamics of swarms viewed as
  complex systems.
\newblock {\em Math. Models Methods Appl. Sci.}, 22(suppl. 1):1140006, 29,
  2012.

\bibitem{ExplicitMPC1}
A.~Bemporad, M.~Morari, V.~Dua, and E.~N. Pistikopoulos.
\newblock The explicit linear quadratic regulator for constrained systems.
\newblock {\em Automatica}, 38(1):3--20, 2002.

\bibitem{LinMPC}
J.~Berberich, J.~Köhler, M.~A. Müller, and F.~Allgöwer.
\newblock Linear tracking {MPC} for nonlinear systems—part i: The model-based
  case.
\newblock {\em IEEE Transactions on Automatic Control}, 67(9):4390--4405, 2022.

\bibitem{blondel2009krause}
V.~D. Blondel, J.~M. Hendrickx, and J.~N. Tsitsiklis.
\newblock On krause's multi-agent consensus model with state-dependent
  connectivity.
\newblock {\em IEEE transactions on Automatic Control}, 54(11):2586--2597,
  2009.

\bibitem{MR4046175}
M.~Burger, R.~Pinnau, C.~Totzeck, O.~Tse, and A.~Roth.
\newblock Instantaneous control of interacting particle systems in the
  mean-field limit.
\newblock {\em J. Comput. Phys.}, 405:109181, 20, 2020.

\bibitem{canizo2011well}
J.~A. Canizo, J.~A. Carrillo, and J.~Rosado.
\newblock A well-posedness theory in measures for some kinetic models of
  collective motion.
\newblock {\em Math. Models Methods Appl. Sci.}, 21(03):515--539, 2011.

\bibitem{Frankowska2018}
P.~Cannarsa and H.~Frankowska.
\newblock Value function, relaxation, and transversality conditions in infinite
  horizon optimal control.
\newblock {\em Journal of Mathematical Analysis and Applications},
  457(2):1188--1217, 2018.
\newblock Special Issue on Convex Analysis and Optimization: New Trends in
  Theory and Applications.

\bibitem{CFPT15}
M.~Caponigro, M.~Fornasier, B.~Piccoli, and E.~Tr\'{e}lat.
\newblock Sparse stabilization and control of alignment models.
\newblock {\em Math. Models Methods Appl. Sci.}, 25(3):521--564, 2015.

\bibitem{carrillo2014derivation}
J.~A. Carrillo, Y.-P. Choi, and M.~Hauray.
\newblock The derivation of swarming models: mean-field limit and {W}asserstein
  distances.
\newblock In {\em Collective dynamics from bacteria to crowds}, pages 1--46.
  Springer, 2014.

\bibitem{carrillo2010particle}
J.~A. Carrillo, M.~Fornasier, G.~Toscani, and F.~Vecil.
\newblock Particle, kinetic, and hydrodynamic models of swarming.
\newblock In {\em Mathematical modeling of collective behavior in
  socio-economic and life sciences}, pages 297--336. Springer, 2010.

\bibitem{carrillo2021controlling}
J.~A. Carrillo, D.~Kalise, F.~Rossi, and E.~Tr\'{e}lat.
\newblock Controlling swarms toward flocks and mills.
\newblock {\em SIAM Journal on Control and Optimization}, 60(3):1863--1891,
  2022.

\bibitem{CKPP19}
Y.-P. Choi, D.~Kalise, J.~Peszek, and A.~A. Peters.
\newblock A collisionless singular {C}ucker-{S}male model with decentralized
  formation control.
\newblock {\em SIAM J. Appl. Dyn. Syst.}, 18(4):1954--1981, 2019.

\bibitem{clarke_vinter}
F.~H. Clarke and R.~B. Vinter.
\newblock The relationship between the maximum principle and dynamic
  programming.
\newblock {\em SIAM Journal on Control and Optimization}, 25(5):1291--1311,
  1987.

\bibitem{MR2165531}
S.~Cordier, L.~Pareschi, and G.~Toscani.
\newblock On a kinetic model for a simple market economy.
\newblock {\em J. Stat. Phys.}, 120(1-2):253--277, 2005.

\bibitem{couzin2005effective}
I.~D. Couzin, J.~Krause, N.~R. Franks, and S.~A. Levin.
\newblock Effective leadership and decision-making in animal groups on the
  move.
\newblock {\em Nature}, 433(7025):513--516, 2005.

\bibitem{MR3308728}
E.~Cristiani, B.~Piccoli, and A.~Tosin.
\newblock {\em Multiscale modeling of pedestrian dynamics}, volume~12 of {\em
  MS\&A. Model. Simul. Appl.}
\newblock Springer, Cham, 2014.

\bibitem{cucker2007emergent}
F.~Cucker and S.~Smale.
\newblock Emergent behavior in flocks.
\newblock {\em IEEE Transactions on automatic control}, 52(5):852--862, 2007.

\bibitem{cucker2004modeling}
F.~Cucker, S.~Smale, and D.-X. Zhou.
\newblock Modeling language evolution.
\newblock {\em Foundations of Computational Mathematics}, 4:315--343, 2004.

\bibitem{MR3119732}
P.~Degond, J.-G. Liu, S.~Motsch, and V.~Panferov.
\newblock Hydrodynamic models of self-organized dynamics: derivation and
  existence theory.
\newblock {\em Methods Appl. Anal.}, 20(2):89--114, 2013.

\bibitem{LucaSDRE}
S.~Dolgov, D.~Kalise, and L.~Saluzzi.
\newblock Optimizing semilinear representations for state-dependent riccati
  equation-based feedback control.
\newblock {\em IFAC-PapersOnLine}, 55(30):510--515, 2022.
\newblock 25th International Symposium on Mathematical Theory of Networks and
  Systems MTNS 2022.

\bibitem{DKS23}
S.~Dolgov, D.~Kalise, and L.~Saluzzi.
\newblock Data-driven tensor train gradient cross approximation for
  hamilton–jacobi–bellman equations.
\newblock {\em SIAM Journal on Scientific Computing}, 45(5):A2153--A2184, 2023.

\bibitem{dyer2009leadership}
J.~R. Dyer, A.~Johansson, D.~Helbing, I.~D. Couzin, and J.~Krause.
\newblock Leadership, consensus decision making and collective behaviour in
  humans.
\newblock {\em Philos. Trans. Roy. Soc. B}, 364(1518):781--789, 2009.

\bibitem{d2006self}
M.~R. D’Orsogna, Y.-L. Chuang, A.~L. Bertozzi, and L.~S. Chayes.
\newblock Self-propelled particles with soft-core interactions: patterns,
  stability, and collapse.
\newblock {\em Phys. Rev. Lett.}, 96(10):104302, 2006.

\bibitem{6315560}
A.~{Eqtami}, D.~V. {Dimarogonas}, and K.~J. {Kyriakopoulos}.
\newblock Event-based model predictive control for the cooperation of
  distributed agents.
\newblock In {\em 2012 Amer. Control Conf.}, pages 6473--6478, 2012.

\bibitem{Giselle}
G.~Estrada-Rodriguez and H.~Gimperlein.
\newblock Interacting particles with {L}\'{e}vy strategies: limits of transport
  equations for swarm robotic systems.
\newblock {\em SIAM J. Appl. Math.}, 80(1):476--498, 2020.

\bibitem{MR4028474}
M.~Fornasier, S.~Lisini, C.~Orrieri, and G.~Savar\'{e}.
\newblock Mean-field optimal control as gamma-limit of finite agent controls.
\newblock {\em European J. Appl. Math.}, 30(6):1153--1186, 2019.

\bibitem{MR3264236}
M.~Fornasier and F.~Solombrino.
\newblock Mean-field optimal control.
\newblock {\em ESAIM Control Optim. Calc. Var.}, 20(4):1123--1152, 2014.

\bibitem{Meurer}
G.~Freudenthaler and T.~Meurer.
\newblock P{DE}-based multi-agent formation control using flatness and
  backstepping: analysis, design and robot experiments.
\newblock {\em Automatica}, 115:108897, 13, 2020.

\bibitem{Garnier}
J.~Garnier, G.~Papanicolaou, and T.-W. Yang.
\newblock Consensus convergence with stochastic effects.
\newblock {\em Vietnam J. Math.}, 45(1-2):51--75, 2017.

\bibitem{MR2887663}
J.~G\'{o}mez-Serrano, C.~Graham, and J.-Y. Le~Boudec.
\newblock The bounded confidence model of opinion dynamics.
\newblock {\em Math. Models Methods Appl. Sci.}, 22(2):1150007, 46, 2012.

\bibitem{HaHaKim}
S.-Y. Ha, T.~Ha, and J.-H. Kim.
\newblock Emergent behavior of a cucker-smale type particle model with
  nonlinear velocity couplings.
\newblock {\em IEEE Transactions on Automatic Control}, 55(7):1679--1683, 2010.

\bibitem{RoM3}
J.~Hahn, U.~Kruger, and T.~F. Edgar.
\newblock Application of model reduction for model predictive control.
\newblock {\em IFAC Proceedings Volumes}, 35(1):393--398, 2002.
\newblock 15th IFAC World Congress.

\bibitem{han2017resolving}
Y.~Han, A.~Hegyi, Y.~Yuan, S.~Hoogendoorn, M.~Papageorgiou, and C.~Roncoli.
\newblock Resolving freeway jam waves by discrete first-order model-based
  predictive control of variable speed limits.
\newblock {\em Transportation Research Part C: Emerging Technologies},
  77:405--420, 2017.

\bibitem{herty2018suboptimal}
M.~Herty and D.~Kalise.
\newblock Suboptimal nonlinear feedback control laws for collective dynamics.
\newblock In {\em 2018 IEEE 14th International Conference on Control and
  Automation (ICCA)}, pages 556--561. IEEE, 2018.

\bibitem{MR2580958}
M.~Herty and L.~Pareschi.
\newblock Fokker-{P}lanck asymptotics for traffic flow models.
\newblock {\em Kinet. Relat. Models}, 3(1):165--179, 2010.

\bibitem{herty2015mean}
M.~Herty, L.~Pareschi, and S.~Steffensen.
\newblock Mean--field control and {R}iccati equations.
\newblock {\em Netw. Heterog. Media}, 10(3):699, 2015.

\bibitem{MR2861587}
M.~Herty and C.~Ringhofer.
\newblock Averaged kinetic models for flows on unstructured networks.
\newblock {\em Kinet. Relat. Models}, 4(4):1081--1096, 2011.

\bibitem{optPDE_Kunisch}
R.~Herzog and K.~Kunisch.
\newblock Algorithms for pde-constrained optimization.
\newblock {\em GAMM-Mitteilungen}, 33(2):163--176, 2010.

\bibitem{LMPC2}
L.~Hewing, K.~P. Wabersich, M.~Menner, and M.~N. Zeilinger.
\newblock Learning-based model predictive control: Toward safe learning in
  control.
\newblock {\em Annual Review of Control, Robotics, and Autonomous Systems},
  3(1):269--296, 2020.

\bibitem{RoM2}
S.~Hovland, K.~Willcox, and J.~T. Gravdahl.
\newblock {MPC} for large-scale systems via model reduction and multiparametric
  quadratic programming.
\newblock In {\em Proceedings of the 45th IEEE Conference on Decision and
  Control}, pages 3418--3423, 2006.

\bibitem{huang2003individual}
M.~Huang, P.~E. Caines, and R.~P. Malham{\'e}.
\newblock Individual and mass behaviour in large population stochastic wireless
  power control problems: centralized and nash equilibrium solutions.
\newblock In {\em 42nd IEEE International Conference on Decision and Control
  (IEEE Cat. No. 03CH37475)}, volume~1, pages 98--103. IEEE, 2003.

\bibitem{AstolfiSDRE}
A.~Jones and A.~Astolfi.
\newblock On the solution of optimal control problems using parameterized
  state-dependent {R}iccati equations.
\newblock In {\em 2020 59th IEEE Conference on Decision and Control (CDC)},
  pages 1098--1103, 2020.

\bibitem{BioKing}
A.~J. King, S.~J. Portugal, D.~Strömbom, R.~P. Mann, J.~A. Carrillo,
  D.~Kalise, G.~de~Croon, H.~Barnett, P.~Scerri, R.~Groß, D.~R. Chadwick, and
  M.~Papadopoulou.
\newblock Biologically inspired herding of animal groups by robots.
\newblock {\em Methods in Ecology and Evolution}, 14(2):478--486, 2023.

\bibitem{Adam}
D.~P. Kingma and J.~Ba.
\newblock Adam: {A} method for stochastic optimization.
\newblock In Y.~Bengio and Y.~LeCun, editors, {\em 3rd International Conference
  on Learning Representations, {ICLR} 2015, San Diego, CA, USA, May 7-9, 2015,
  Conference Track Proceedings}, 2015.

\bibitem{GP2}
J.~Kocijan and R.~Murray-Smith.
\newblock {\em Nonlinear Predictive Control with a Gaussian Process Model},
  pages 185--200.
\newblock Springer Berlin Heidelberg, Berlin, Heidelberg, 2005.

\bibitem{GP1}
J.~Kocijan, R.~Murray-Smith, C.~Rasmussen, and B.~Likar.
\newblock Predictive control with gaussian process models.
\newblock In {\em The IEEE Region 8 EUROCON 2003. Computer as a Tool.},
  volume~1, pages 352--356 vol.1, 2003.

\bibitem{RL1}
S.~Lale, K.~Azizzadenesheli, B.~Hassibi, and A.~Anandkumar.
\newblock Model learning predictive control in nonlinear dynamical systems.
\newblock In {\em 2021 60th IEEE Conference on Decision and Control (CDC)},
  pages 757--762, 2021.

\bibitem{RL2}
D.~Limon, J.~Calliess, and J.~Maciejowski.
\newblock Learning-based nonlinear model predictive control.
\newblock {\em IFAC-PapersOnLine}, 50(1):7769--7776, 2017.
\newblock 20th IFAC World Congress.

\bibitem{relax_DP}
B.~Lincoln and A.~Rantzer.
\newblock Relaxing dynamic programming.
\newblock {\em IEEE Transactions on Automatic Control}, 51(8):1249--1260, 2006.

\bibitem{offline2}
E.~Maddalena, C.~{da S. Moraes}, G.~Waltrich, and C.~Jones.
\newblock A neural network architecture to learn explicit {MPC} controllers
  from data.
\newblock {\em IFAC-PapersOnLine}, 53(2):11362--11367, 2020.
\newblock 21st IFAC World Congress.

\bibitem{Gp3}
M.~Maiworm, D.~Limon, and R.~Findeisen.
\newblock Online learning-based model predictive control with gaussian process
  models and stability guarantees.
\newblock {\em International Journal of Robust and Nonlinear Control},
  31(18):8785--8812, 2021.

\bibitem{RL3}
J.~E. Morinelly and B.~E. Ydstie.
\newblock Dual {MPC} with reinforcement learning.
\newblock {\em IFAC-PapersOnLine}, 49(7):266--271, 2016.
\newblock 11th IFAC Symposium on Dynamics and Control of Process
  SystemsIncluding Biosystems DYCOPS-CAB 2016.

\bibitem{motsch2014heterophilious}
S.~Motsch and E.~Tadmor.
\newblock Heterophilious dynamics enhances consensus.
\newblock {\em SIAM review}, 56(4):577--621, 2014.

\bibitem{norouzi2023integrating}
A.~Norouzi, H.~Heidarifar, H.~Borhan, M.~Shahbakhti, and C.~R. Koch.
\newblock Integrating machine learning and model predictive control for
  automotive applications: A review and future directions.
\newblock {\em Engineering Applications of Artificial Intelligence},
  120:105878, 2023.

\bibitem{KPAsurvey15}
K.-K. Oh, M.-C. Park, and H.-S. Ahn.
\newblock A survey of multi-agent formation control.
\newblock {\em Automatica}, 53:424--440, 2015.

\bibitem{MR3157726}
A.~A. Peters, R.~H. Middleton, and O.~Mason.
\newblock Leader tracking in homogeneous vehicle platoons with broadcast
  delays.
\newblock {\em Automatica}, 50(1):64--74, 2014.

\bibitem{offline1}
S.~S. {Pon Kumar}, A.~Tulsyan, B.~Gopaluni, and P.~Loewen.
\newblock A deep learning architecture for predictive control.
\newblock {\em IFAC-PapersOnLine}, 51(18):512--517, 2018.
\newblock 10th IFAC Symposium on Advanced Control of Chemical Processes ADCHEM
  2018.

\bibitem{pontryagin1962}
L.~S. Pontryagin, V.~G. Boltyanskii, R.~V. Gamkrelidze, and E.~F. Mishchenko.
\newblock {\em The mathematical theory of optimal processes, Chapter 1}.
\newblock Classics of Soviet Mathematics. Interscience Publishers John Wiley \&
  Sons, Inc., 1962.
\newblock Translated from the Russian by K. N. Trirogoff; edited by L. W.
  Neustadt.

\bibitem{stern2018dissipation}
R.~E. Stern, S.~Cui, M.~L. Delle~Monache, R.~Bhadani, M.~Bunting, M.~Churchill,
  N.~Hamilton, H.~Pohlmann, F.~Wu, B.~Piccoli, et~al.
\newblock Dissipation of stop-and-go waves via control of autonomous vehicles:
  Field experiments.
\newblock {\em Transp. Research Part C: Emerging Techn.}, 89:205--221, 2018.

\bibitem{subbotina2006method}
N.~Subbotina.
\newblock The method of characteristics for {H}amilton—{J}acobi equations and
  applications to dynamical optimization.
\newblock {\em Journal of mathematical sciences}, 135:2955--3091, 2006.

\bibitem{MR2247927}
G.~Toscani.
\newblock Kinetic models of opinion formation.
\newblock {\em Commun. Math. Sci.}, 4(3):481--496, 2006.

\bibitem{MR3948232}
A.~Tosin and M.~Zanella.
\newblock Kinetic-controlled hydrodynamics for traffic models with
  driver-assist vehicles.
\newblock {\em Multiscale Model. Simul.}, 17(2):716--749, 2019.

\bibitem{FastMPC}
I.~J. Wolf and W.~Marquardt.
\newblock Fast {NMPC} schemes for regulatory and economic {NMPC} – a review.
\newblock {\em Journal of Process Control}, 44:162--183, 2016.

\bibitem{RoM1}
T.~Zhao, Y.~Zheng, J.~Gong, and Z.~Wu.
\newblock Machine learning-based reduced-order modeling and predictive control
  of nonlinear processes.
\newblock {\em Chemical Engineering Research and Design}, 179:435--451, 2022.

\bibitem{zhou1990maximum}
X.~Zhou.
\newblock Maximum principle, dynamic programming, and their connection in
  deterministic control.
\newblock {\em Journal of Optimization Theory and Applications},
  65(2):363--373, 1990.

\end{thebibliography}
	
\end{document}